\documentclass[11pt]{article}

  \usepackage[utf8]{inputenc}
  \usepackage{amssymb,amsthm,amsmath}
  \usepackage[hidelinks]{hyperref}
  \usepackage{float}
  \usepackage{caption}
  \usepackage{subfig}
  \usepackage{graphicx}
  \usepackage{multirow}

\theoremstyle{plain}
  \newtheorem{theorem}{Theorem}  
  \newtheorem{lemma}{Lemma}

\allowdisplaybreaks[2]

\begin{document}

\title{Nonparametric Multiple Change Point Detection \\for Non-Stationary Time Series}
\author{Zixiang Guan, Gemai Chen  \\
Department of Mathematics and Statitsics \\University of Calgary, Calgary, Alberta, T2N 1N4 \\
Email: \texttt{zixiang.guan@ucalgary.ca} \\ \texttt{cheg@ucalgary.ca} }

\maketitle

\newpage

\mbox{}
\vspace*{2in}
\begin{center}
\textbf{Author's Footnote:}
\end{center}
Zixiang Guan is PhD, Department of Mathematics and Statistics, University of Calgary (Email: \texttt{zixiang.guan@ucalgary.ca}),
Gemai Chen is Professor (Email: \texttt{gchen@math.ucalgary.ca}), Department of Mathematics and Statistics, University of Calgary, Alberta, Canada.

\newpage
\begin{center}
\textbf{Abstract}
\end{center}
This article considers a nonparametric method for detecting change points in non-stationary time series. The proposed method will divide the time series into several segments so that between two adjacent segments, the normalized spectral density functions are different. The theory is based on the assumption that within each segment, time series is a linear process, which means that our method works not only for classic time series models, e.g., causal and invertible ARMA process, but also preserves good performance for non-invertible moving average process. We show that our estimations for change points are consistent. Also, a Bayesian information criterion is applied to estimate the member of change points consistently. Simulation results as well as empirical results will be presented.

\vspace*{.3in}


\newpage

\section{Introduction}
Time series analysis is a well-developed branch of statistics with a wide range of applications in engineering, economics, biology and so on. Generally speaking, when investigating the theoretical properties as well as analyzing real data both in time domain and frequency domain, it is often assumed that time series is stationary. However in application, stationarity may be violated. How to analyze non-stationary time series is challenging. In ARIMA model ( Brockwell and Davis 1991), data is differenced finite times so that it reduces to ARMA process. In another way, we may assume that non-stationary process consists of several stationary ones. Our goal is to segment time series properly.\par

 Consider a sequence of data $\{X_1,\ldots,X_N\}$ and let $\tau^0_0,\tau^0_1,\tau_2^0,\ldots,\tau_K^0,\tau^0_{K+1}$ be nonnegative integers satisfying $0=\tau^0_0<\tau^0_1<\tau_2^0<\ldots<\tau_K^0<\tau^0_{K+1}=N$. We assume that within the $j$th segment of data, i.e., $\tau_{j-1}^0+1\leq t\leq \tau_j^0$, $1\leq j\leq K+1$, $X_t$ is stationary. $\tau_j^0$ are called structural breaks, or change points, which are unknown. $K$ is the number of change points.  From the perspective of time domain, we can assume that each stationary segment can be modeled by appropriate statistical models while its structure varies across different segments. Davis, Lee and
Rodriguez-Yam (2006) divided non-stationary time series into several different autoregressive processes. Kitagawa and Akaike (1978) detected change points by AIC criterion. In frequency domain, Ombao, Raz, Von Sachs and Malow (2001) used a family of orthogonal wavelet called SLEX to partition non-stationary process into stationary ones.  Lavielle and Lude\~{n}a (2000) estimated change points using Whittle log-likelihood when time series was parametric. In Korkas and Fryzlewicz (2017), Locally Stationary Wavelets was applied to estimate the second-order structure, then Wild Binary Segmentation was imposed to divide the time series into several segments based on CUSUM statistics. \par

When autocovariance function is absolutely summable, it is well known that stationary process has spectral density function, which is the Fourier transformation applied to autocovariance function (Brockwell and Davis 1991). Spectral density function preserves a good property which is that after being properly normalized, it becomes a well-defined probability density function (Priestley 1981). In probability theory, there are some existing functions to measure the difference between probability density functions. Kullback-Leibler divergence (K-L divergence) is one of them. \par
In  Kullback and Leibler (1951), a divergence function was introduced to measure the discrimination information between two distribution functions. The original definition is as follows. Suppose there are two probability distributions, $Z_1(x)$ and $Z_2(x)$. $z_1(x)$ and $z_2(x)$ denote the probability density functions of $Z_1$ and $Z_2$, respectively. Then Kullback-Leibler (K-L) divergence from $Z_2$ to $Z_1$ is
\begin{center}
$\displaystyle D_{\mathrm{KL}}(Z_1\|Z_2)=\int z_1(x)\log\left(\frac{z_1(x)}{z_2(x)}\right)dx$.
\end{center}
As we can see, K-L divergence is defined for probability density function, which is non-negative. Since spectral density function is also non-negative, we will generalize K-L divergence so that it is applicable to non-negative functions, which is still called Kullback-Leibler divergence in this article and will be applied later in the forthcoming sections. The definition is as follows. Suppose we have two non-negative functions, $f_1$ and $f_2$, defined on a common support. Then the Kullback-Leibler divergence of $f_1(x)$ with respect to $f_2(x)$ is
\begin{center}
$\displaystyle D_{\mathrm{KL}}(f_1\|f_2)=\int f_1(x)\log\frac{st(f_1)}{st(f_2)}dx$,
\end{center}
where $st(f_1)=\frac{f_1(x)}{F_1}$, $st(f_2)=\frac{f_2(x)}{F_2}$, and $F_1=\int f_1(x)dx$, $F_2=\int f_2(x)dx$.
By Gibbs's inequality, this new K-L divergence is still non-negative and is equal to 0 if and only if $f_1=cf_2$ almost everywhere, where $c$ is an appropriate constant. In this article, we apply K-L divergence to normalized spectral density functions to define our objective function. Then we estimate change points by maximizing the objective function. That is, we find the locations where discrepancy between different spectral density functions reaches its maximum. When estimating spectral density function, we adopt the classical method. That is, we first calculate periodogram then smooth it by choosing appropriate spectral window. Although we assume that each stationary segment is a linear process, which is a class of stationary time series more general than autoregressive process and ARMA model, we calculate spectral density function without estimating any parameters. So our change point detection method is nonparametric, which we call it nonparametric spectral change-point detection (NSCD). In application, we do not know the number of change points, so a BIC criterion,which is similar to Yao (1988) and Zou, Yin, Feng, and Wang (2014), is proposed. The consistency of both change point estimation and BIC criterion can be guaranteed. Dynamic programming algorithm (Hawkins 2001) is used, but due to its computational complexity, we adopt the screening algorithm (Zou \textsl{et al}. 2014). Also, Pruned Exact Linear Time (PELT) by Killick, Fearnhead and
Eckley (2012) can also boost the speed of algorithm when estimating the number and locations of change points simultaneously.\par

The rest of this article is organized as follows. In Section 2, we describe our objective functions. Asymptotic properties are presented in Section 3. Implementation and further details of our algorithm will be given in Section 4. In Section 5, numerical simulations as well as comparison with several methodologies are shown. The analysis of EEG data is presented in Section 6. \par

\section{Model and Methodology \label{ch3s1}}
Consider non-stationary time series $X_t=\sum\limits_{j=-\infty}^{+\infty} a_j(k)\xi_{t-j}$, $\xi_j\overset{i.i.d}{\sim}(0,\sigma^2)$, $\tau^0_{k-1}< t\leq\tau^0_k$, $k=1,\ldots,K+1$ with $\tau_0^0=0$, $\tau^0_{K+1}=N$, $\tau_k^0-\tau_{k-1}^0=N_k$. Here we assume that all $\xi_t$ are independent and identically distributed. $\{a_j(k)\}$ satisfy $\sum\limits_{j}|a_j(k)|<+\infty$, $\forall k$, so that when $\tau_{k-1}^0< t\leq\tau_k^0$, spectral density function exists, denoted by $f_k$ (Brockwell and Davis 1991). If $EX_t\neq0$, we can subtract mean from observations by $X_t-\mu$. Otherwise, we assume that $EX_t=0$ for all $t$. Our goal is to detect $\tau_k$ when $f_k$ and $f_{k+1}$ are different. If we could find a spectral density function $f$ that it is different from all $f_k$, then a discriminant function can be applied so that it reaches its maximum at $\tau_k$. Here we adopt Kullback-Leibler divergence and define our objective function as below:
\begin{center}
$\displaystyle R(\tau_1,\ldots,\tau_{K+1})=\sum\limits_{k=1}^{K+1}(\tau_k-\tau_{k-1})\int_{-\pi}^\pi \hat{f}_k(u)\log\frac{st(\hat{f}_k)}{st(f)}du$.
\end{center}
Here $\hat{F}_k=\int_{-\pi}^\pi \hat{f}_k(u)du$, $F=\int_{-\pi}^\pi f(u)du$, $st(\hat{f}_k) = \frac{\hat{f}_k}{\hat{F}_k}$, and $st(f)=\frac{f}{F}$ so that $st(\hat{f}_k)$ and $st(f)$ become probability density functions, and $0<\tau_1<\tau_2<\cdots<\tau_K<N$ is a possible partition of the time series.
There are two ways to find $f$. We can give an estimation of the spectral density function $\hat{f}_0$, based on the whole time series. Although the whole time series is non-stationary, the estimated spectral density function still converges to a well-defined function, which is the weighted sum of $f_k$ (See Lemma \ref{ch3l2} for details). If we know that data is not white noise, then $st(f)=\frac{1}{2\pi}$. \par

There are literatures concerning the application of K-L divergence in time series. Parzen (1982) used it to estimate the parameters of autoregressive process. Parzen (1983) extended it to estimate ARMA process. Shore (1981) calculated the minimum K-L divergence to estimate spectrum given its priori spectrum has an exponential form.  See Rao (1993) for more reviews. From information point of view, K-L divergence measures the information gain when a new probability density function is used instead of the old one. So our objective function will detect the locations where we maximize the information gain when a new spectral density function $f$ comes in. \par
To estimate the spectral density function, the classical methodology is adopted here. First, periodogram will be calculated as follows:
\begin{center}
$\displaystyle I_{\tau_{k-1}+1,\tau_k}(\lambda)=\frac{1}{\tau_k-\tau_{k-1}}\left|\sum\limits_{j=\tau_{k-1}+1}^{\tau_k}X_je^{-ij\lambda}\right|^2$
\end{center}
Since periodogram is not consistent, we choose a spectral window, $W(u)$, to smooth periodogram, then
\begin{center}
$\displaystyle \hat{f}_k(\lambda)=m\int_{-\pi}^{\pi} W(m(\lambda-u))I_{\tau_{k-1}+1,\tau_k}(u)du$, $\lambda\in[-\pi,\pi]$.
\end{center}
There are several choices for $W(u)$ (Priestley 1981). Since we want the normalized $\hat{f}_k$ to be a probability density function, i.e., $\hat{f}_k\geq 0$, Bartlett kernel is chosen, which is $mW(mu)=\sin^2(mu/2)/(2\pi m \sin^2(u/2))$, where $m$ is the bandwidth.
One usually applies Fourier transformation to achieve an estimator on a grid of frequencies, denoted by $\Lambda=\{\lambda_1,\ldots,\lambda_{N_\lambda}\}$, where $N_\lambda$ is the cardinality of $\Lambda$, which could tend to infinity as $N$ goes to infinity. Here we still use integral to denote the summation of estimators. We may set $N_\lambda=N$ so $N_\lambda$ could tend to infinity, or $\Lambda$ could be Nyquist frequency. Based on this grid of frequency, K-L divergence applied to normalized spectral density function is still non-negative, and equals zero if and only if two normalized spectral density functions are equal on $[-\pi,\pi]$.\par
When estimating change points, we have no idea about the number of change points in the data, so $K$ should be estimated. Yao and Au (1989), Zou \textsl{et al}. (2014) gave BIC criterion and showed its consistency. Following their work, we also propose a BIC criterion as follows:
\begin{center}
$BIC_L=-\max_{\tau_1^\prime,\ldots,\tau_K^\prime}+LC_N$.
\end{center}
Here $C_N$ is an appropriate constant which will be illustrated later. So our estimator $\hat{K}$ is chosen by minimizing the criterion above.

\section{Asymptotic Theory \label{ch3s2}}
In Section \ref{ch3s1}, we can see that $K$ is a constant, and $\tau_i^0$, $i=1,2,\ldots, K$, change with $N$. So there is a sequence of constants, $0<\kappa_1^0<\kappa_2^0<\cdots<\kappa_K^0<1$, such that $\{X_t,t=1,2,\ldots, 1\}$ is a realization of non-stationary time series with $\tau_k^0=[\kappa_k^0N]$, $k=1,2,\ldots, K$, where $[x]$ denotes the largest integer which is not greater than $x$. Their estimators are denoted by $\hat{\kappa}_k$, $k=1,\ldots,K$. We can estimate $\tau_k^0$ first, denoted by $\hat{\tau}_k$, by maximizing the objective function $R(\tau_1,\ldots,\tau_K)$, then $\hat{\kappa}_k=\hat{\tau}_k/N$. In literature, Yao and Au (1989) achieved a consistent estimation of $O_p(1)$ for $\hat{\tau}_k$, which means that the difference between estimators and true change points is no bigger than a constant.  Zou \textsl{et al}. (2014) also drew the same conclusion when the number of change points was constant. In  Davis \textsl{et al}. (2006), the consistency was attained in the sense that $|\hat{\kappa}_k-\kappa_k|<\epsilon$ in probability 1, where $\epsilon$ was some constant. The reason that estimators cannot converge as fast as those in Yao and Au (1989) and Zou \textsl{et al}. (2014) is that estimating AR process needs a sufficient number of samples. To guarantee estimate accuracy, we always find the next change point which is $ml$ away from the previous one. For example, after giving $\hat{\tau}_j$, we look for the next change point starting from $\hat{\tau}_j+ml$. Our results are similar to Davis \textsl{et al}. (2006) and based on a set of discrete frequencies $\Lambda$. The following assumptions are needed to obtain the consistency:
\begin{itemize}
\item[A1:] $E\xi_t^{8q}<\infty$ where $q$ is some integer satisfying $q\geq3$, $\{a_j(k)\}$ converge absolutely, $\forall k$.
\item[A2:] $W(v)$ is a non-negative, even, bounded, integrable function, $\int_{-\pi}^{\pi}W(v)dv=1$,\\ $\int_{-\pi}^\pi \left(W(v)\right)^{1-\frac{1}{2q}}dv<\infty$.
\item[A3:] $N_{\mathrm{min}}/N\rightarrow c_{\mathrm{min}}>0$, as $N\rightarrow\infty$, and $N_{\mathrm{min}}>m$, where $N_{\mathrm{min}}=\min_{1\leq k\leq K}(\tau_k^0-\tau^0_{k-1})$.
\item[A4:] $\forall k$, $f_k$ is everywhere positive and satisfies uniform Lipschitz condition:
\begin{center}
$\left|f_k(u_1)-f_k(u_2)\right|\leq B_{f_k}|u_1-u_2|$
\end{center}
$\forall u_1\in[-\pi,\pi]$, $u_2\in[-\pi,\pi]$, where $B_{f_k}$ is a constant.
\item[A5:] $m=O(N^\alpha)$, where $\frac{1}{4}\leq\alpha<\alpha+\frac{3}{2q}<\frac{1}{2}$.
\item[A6:] $w(0)=1$, and $w(v)$ has continuous derivatives to the order of $2q$.
\item[A7:] $\{f_k(u)\}$ are linearly independent for $u \in[-\pi,\pi]$.
\item[A8:] $ml=c_{\mathrm{ml}}N$, where $c_{\mathrm{ml}}>0$. $ml<N_{\mathrm{min}}$
\end{itemize}
In A8, $ml$ is the minimal length of time series when estimating spectral density functions, since a sufficient number of observations is always necessary, especially when the convergence rate of estimators is slow. Also, $ml<N_{min}$ so that all change points are distinguishable. Assumption 7 is given to guarantee that any linear combination is not equal to $f_j$, $\forall j=1,\ldots,K$.
Assumption 1-6 are similar to those in Woodroofe and Van Ness (1967) so that the $\max\limits_{\lambda_j\in\Lambda}\frac{\hat{f}_k(\lambda_j)}{f_k(\lambda_j)}$ can be bounded in probability. Assumption 4 can be stronger so that in Assumption 5, $\alpha$ can be less than $\frac{1}{4}$ (see Woodroofe and Van Ness (1967) for further details). Theorem \ref{esticonsis3} gives the consistency of change point estimation.
\begin{theorem}\label{esticonsis3}
When $K$ is known, $\hat{\kappa}_j\overset{p}{\rightarrow}\kappa_j$, $\forall j=1,\ldots, K$.
\end{theorem}
To estimate the number of change points, a pre-specified upper bound $K_{\mathrm{max}}$ satisfying $K<K_{\mathrm{max}}$ will be given. Then BIC values for each $1\leq L\leq K_{\mathrm{max}}$ are calculated and $\hat{K}$ will be the number where BIC reaches its minimum. The following theorem establishes the consistency of estimation of BIC criterion.
\begin{theorem}\label{BICconsist3}
If $C_N/N\rightarrow0$, $C_N/N^{\frac{q+3+2q\alpha}{2q}}\rightarrow\infty$, we have $P(\hat{K}=K)\rightarrow1$.
\end{theorem}

\section{Algorithm}

Similar to  Zou \textsl{et al}. (2014), our objective function is separable. For change point detection, a commonly adopted algorithm, Dynamic Programming (Hawkins 2001), can be applied. The main idea of Dynamic Programming is that estimation of $\hat{\tau}_K$ is computed first, which is the rightmost change point. Then the time series data from $1$ to $\hat{\tau}_{K}$ will be divided into $K-1$ parts, and we will estimate $\tau_{K-1}$ recursively. However, the computation complexity is $O(KN^2)$, and taking Discrete Fourier transformation and spectrum smoothing into consideration, it is time-consuming. \par
To reduce computational complexity, Zou \textsl{et al}. (2014) proposed a screening algorithm. For $X_j,\ldots,X_{j+l}$, where $l$ is some constant integer, calculate the location where the function below reaches its maximum.
\begin{eqnarray*}
r_j&=&\arg\max\limits_rr\int_{-\pi}^\pi \hat{f}_{j,j+r}(u)\log\frac{st(\hat{f}_{j,j+r})}{st(f)}du\\
& &+(l-r)\int_{-\pi}^\pi \hat{f}_{j+r+1,j+l}\log\frac{st(\hat{f}_{j+r+1,j+l})}{st(f)}du.
\end{eqnarray*}
Here $\hat{f}_{j,j+r}$ and $\hat{f}_{j+r+1,j+l}$ denote the estimated spectral density functions based on observations from $j$ to $j+r$, and $j+r+1$ to $j+l$. Let $j$ change from $1$ to $n-l$, then we have a set $A_{sc}$ containing all $r_j$, then apply Dynamic Programming on $A_{sc}$. The main idea is that if a change point is included in $X_j,\ldots,X_{j+l}$, then the equation above should reach its maximum at this true change point. That is, $A_{sc}$ contains true change points. Here we should choose $l<N_{\mathrm{min}}$ so that $X_j,\ldots,X_{j+l}$ contain only one change point. \par

Here $\hat{f}_{j,j+r}$ and $\hat{f}_{j+r+1,j+l}$ denote the estimated spectral density function of samples from $j$ to $j+r$, and $j+r+1$ to $j+l$. Let $j$ change from $1$ to $n-l$, then we have a set $A_{sc}$ containing all $r_j$, then apply dynamic programming on $A_{sc}$. The main idea is that if a change point is included in $X_j,\ldots,X_{j+l}$, then the equation above should reach its maximum at this true change point. That is, $A_{sc}$ contains true change points. Here we should choose $l<N_{\mathrm{min}}$ so that $X_j,\ldots,X_{j+l}$ contain only one change point. \par

When calculating the spectrums within each segment of time series, the most widely used method is Fast Fourier Transformation (FFT). However, in FFT, a problem is that spectral density function is estimated on Nyguist frequency. If so, time series with different length is estimated on different set of Fourier frequency. So when calculating the integral in $R(\tau_1,\ldots,\tau_K)$, we cannot align the frequencies where spectral density functions $f_k$ and $f$ are estimated. Our method is that we first choose a set of frequencies, denoted by $\Lambda$, then apply Discrete Fourier Transformation on $\Lambda$ for every subset of samples, which will solve the alignment problem naturally.\par

For BIC criterion, usually Dynamic Programming will be applied first, then the values of BIC criterion for all $L\leq K_{\mathrm{max}}$ will be calculated. Obviously, this will increase complexity. Killick \textsl{et al}. (2012) proposed a method called Pruned Exact Linear Time (PELT), which would significantly reduce computational complexity. The main idea is that when a new sample is included, check all the remaining locations before new sample. If the objective function decreases to the extent that it is larger than the penalty term $C_N$, those locations which do not satisfy the condition will be removed and the next sample will be added into our calculation until the end. The cardinality of the set of all remaining locations is the estimated number of change points and the elements within will be the estimated change points. Under some assumptions, the computational complexity is linear with respect to sample size. \par

In BIC criterion, another problem is how to choose $C_N$. Although we can set $C_N$ to satisfy conditions in Theorem \ref{BICconsist3}, this choice may be too large which leads to underestimation of $K$. To overcome this difficulty, we first choose a length, which equals $ml$, then compute the median, denoted by $me_{\mathrm{BIC}}$, for all values of the function below:
\begin{center}
$\displaystyle \int \hat{f}_{j,j+ml}(u)\log\frac{st(\hat{f}_{j,j+ml})}{st(f)}du$, $\forall j=1,\ldots,K$.
\end{center}
$\hat{f}_{j,j+ml}$ is the estimated spectral density function from $X_j,\ldots,X_{j+ml}$, $\hat{F}_{j,j+ml}$ is its integral. Finally, $C_n=me_{\mathrm{BIC}}\times n^{c}$. Our simulation shows that an appropriate choice for $c$ is 0.73 regardless of $m$. \par
The selection of spectral windows is an important topic in spectral estimation. In Priestley
(1981), bandwidth is selected as follows
\begin{center}
$B_W=2\sqrt{6}\left(\frac{1}{m^r}k^{(r)}\right)^{1/r}$,
\end{center}
where $k^{(r)}=\lim\limits_{u\rightarrow0}\frac{1-w(u)}{|u|^r}$, where $w(u)$ is the inverse Fourier transformation of $W(u)$, and $r$ is the largest integer so that the limit aforementioned exists and is non-zero. $m$ is the scale parameter in spectral windows.  By the proofs of Theorem \ref{esticonsis3}, we can see that estimators of change points reach consistency because $N$ dominate convergence rate. So bandwidth selection does not matter too much, which is verified by our simulations.\par

In application, sometimes sample size is large. Although we can apply some methods, such as screening, PELT, to boost the calculation, the computational complexity is still intolerable. Here, we set a change point searching unit, denoted by $n_{\mathrm{su}}$ (Hawkins 2001). That is, when searching for change points, we add a unit of observations into our calculation each time, not just one observation. This unit is different from $ml$, which is used to calculate spectrums since estimating spectral density function usually needs sufficient amount of observations. By setting this unit, change point can only be estimated at $n_{\mathrm{su}}$, $2n_{\mathrm{su}}$, $3n_{\mathrm{su}}$, and so on. If we set this unit equal to 1, the algorithm degenerates to the general scenario when no searching unit is given. This will dramatically increase the speed of our algorithm. The computational complexity of Dynamic Programming is $O(n^2)$ (Zou \textsl{et al}. 2014), so by setting a searching unit $n_{\mathrm{su}}$, the complexity will drop to $O(n^2/n^2_{\mathrm{su}})$. Apparently the estimation accuracy will be sacrificed, since the true change points and the maximizers of objective function, may not lie on the grid. We suggest the choice of $n_{\mathrm{su}}$ by choosing the desired estimation accuracy first. Intuitively speaking, if we want the estimates having the accuracy of $1\%$ of the total sample size, then we can set $n_{\mathrm{su}}=0.01n$. \par

\section{Simulation \label{ch3s4}}
In this section we show the finite sample properties of our method and compare it with AutoPARM, Wild Binary Segmentation (WBS), Binary Segmentation (BS), MuBred, and NMCD under several cases.


Following Zou \textsl{et al}. (2014),  we calculate the distance between two sets $G$ and $\hat{G}$, which are the true change point set and the estimated set, respectively, by
\begin{center}
$\varrho(\hat{G}||G)=\sup\limits_{b\in G}\inf\limits_{a\in \hat{G}}|a-b|$, and $\varrho(G||\hat{G})=\sup\limits_{a\in \hat{G}}\inf\limits_{b\in G}|a-b|$.
\end{center}
The first measurement shows if there is an estimator close enough to a true change point, while the second measurement reveals the distance between estimates and true change points. When $K$ is known, both measures should give good performances in the sense that $\varrho(\hat{G}||G)$ and $\varrho(G||\hat{G})$ are small. For all the tables, $\varrho(\hat{G}||G)$ and $\varrho(G||\hat{G})$ are shown outside the parentheses while $\varrho(\hat{G}||G)/N$ and $\varrho(G||\hat{G})/N$ are shown inside, which can measure the estimation accuracy of $\hat{\kappa}$. In the following subsections, $\xi_t\overset{\mathrm{i.i.d}}{\sim} N(0,1)$, if it is not mentioned. To guarantee the estimation accuracy of spectral density function, we should set the minimal length of a segment, which is denoted by $ml$. In this simulation, we set $ml=350$ if it is not mentioned. By setting $ml$, we also assume that the distance of two adjacent change points will not be smaller than $ml$, so we set $K_{\mathrm{max}}=6$ for each case. If $K$ should be estimated, we report the percentage when $K$ is accurately detected. All simulations are obtained with 1000 replications. \par
\subsection{AutoRegresstive Process \label{ch3simu1}}
Following the examples in Davis {\sl et al.} (2006), we generate the non-stationary time series from the following. \par
\begin{figure}[H]
\centering
\caption{Normalized Spectral Density Functions in Davis {\sl et al.} (2006)}
\label{Davis}
\includegraphics[scale=0.5]{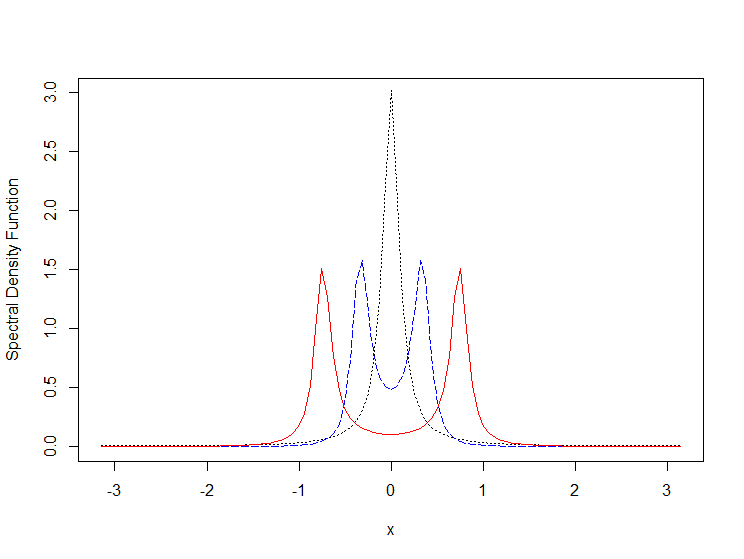}
\caption*{{\small AR(1) is in black dotted line, the second AR is in blue dashed line, the last AR is in red solid line.}}
\end{figure}

Autoregressive Processes (Case 1):
\begin{itemize}
\item[1] $X_t-0.9X_{t-1}=\xi_t$, $1\leq t\leq 1024$,
\item[2] $X_t-1.69X_{t-1}+0.81X_{t-2}=\xi_t$, $1025\leq t\leq 1536$,
\item[3] $X_t-1.32X_{t-1}+0.81X_{t-2}=\xi_t$, $1537\leq t\leq 2048$.
\end{itemize}
and their normalized spectral density functions are shown in Figure \ref{Davis}.
In Table \ref{Davissimu1}, we show the results of our method when $K$ is known with the spectral density function of total samples and white noise as the baseline functions, and different choices of bandwidths, respectively. In Table \ref{Davissimu2}, the simulation is conducted when $K$ is unknown. We adopt Bartlett window since it guarantees the non-negativity of estimated spectral density function. To reduce computational complexity, screening algorithm is applied. \par
\begin{table}[H]
\centering
\caption{Performances of NSCD for Case 1 with Different Baseline Functions and Bandwidths When $K$ Is Known \label{Davissimu1}}
\begin{tabular}{|c|c|c|c|}
\hline
 & & $\varrho(\hat{G}||G)$ & $\varrho(G||\hat{G})$\\
\hline
\multirow{2}{*}{$m=N^{1/3}$} & $\hat{f}_0$ & 22.99 (0.011) & 22.99 (0.011)\\
\cline{2-4}
            & $\frac{1}{2\pi}$ & 23.55 (0.012) & 23.55 (0.012)\\
\hline
\multirow{2}{*}{$m=N^{1/4}$} & $\hat{f}_0$ & 22.61 (0.011) & 22.61 (0.011)\\
\cline{2-4}
            & $\frac{1}{2\pi}$ & 26.47 (0.013) & 26.47 (0.013)\\
\hline
\end{tabular}
\end{table}
\begin{table}[H]
\centering
\caption{BIC criterion of NSCD for AR Processes with Comparison to AutoPARM, WBS, BS, MuBred, and NMCD \label{Davissimu2}}
\begin{tabular}{|c|c|c|c|}
\hline
& $\hat{K}$ & $\varrho(\hat{G}||G)$ & $\varrho(G||\hat{G})$\\
\hline
AutoPARM & 93.3\% & 7.61 (0.004) & 8.05 (0.004)\\
\hline
NSCD & 98.8\% & 25.00 (0.012) & 37.33 (0.018)\\
\hline
MuBred & 97.3 \% & 12.219 (0.006) & 16.229 (0.008)\\
\hline
WBS & 37.9\% & 66.65 (0.033) & 183.184 (0.089)\\
\hline
BS & 62.7\% & 84.875 (0.041) & 122.058 (0.060)\\
\hline
NMCD & 0\% & 147.4365 (0.072) & 554.9088 (0.271)\\
\hline
\end{tabular}
\end{table}
From Table \ref{Davissimu1}, the choice of baseline function does not make much difference. Also, under different bandwidths, results are quite similar. We may suggest to choose larger bandwidth for large sample size since the bias of estimation is small.
From Table \ref{Davissimu2}, it is not surprising that AutoPARM outperforms NSCD, because AutoPARM knows the exact structure of data. WBS and BS perform weaker than NSCD, and cannot estimate the number of change points accurately. The performance of MuBred is better than NSCD, although it does not need the assumption of Autoregressive process. NMCD fails to capture the true change points, since it is designed for independent random variables. In fact, it prefers to overestimate the number of change points, which is 3.86.  \par
\subsection{ARMA Process and Invertible Moving Average Process \label{ch3simu2}}
Next, we simulate from ARMA and invertible MA processes. Table \ref{ARMA1} and Table \ref{ARMA2} contain the results. \par
ARMA (Case 2):
\begin{itemize}
\item[1:] $X_t-X_{t-1}+0.25X_{t-2}=\xi_t+0.8\xi_{t-1}$, $1\leq t\leq 500$,
\item[2:] $X_t-0.5X_{t-1}=\xi_t$, $501\leq t\leq 1100$,
\item[3:] $X_t-1.7X_{t-1}+0.9X_{t-2}-0.168X_{t-3}=\xi_t-1.6\xi_{t-1}+0.79\xi_{t-2}-0.12\xi_{t-3}$, $1101\leq t\leq 1800$.
\end{itemize}
Invertible MA (Case 3):
\begin{itemize}
\item[1:] $X_t=(3+B)(2-B)\xi_t$, $1\leq t\leq 500$,
\item[2:] $X_t=(3-B)(2-B)\xi_t$, $501\leq t\leq 1100$,
\item[3:] $X_t=(3+B)(2-B)\xi_t$, $1101\leq t\leq 1800$.
\end{itemize}
\begin{table}[H]
\centering
\caption{Performance of NSCD for ARMA and Invertible MA Processes with Different Baseline Functions and Bandwidths When $K$ Is Known \label{ARMA1}}
\begin{tabular}{|c|c|c|c|c|c|}
\hline
 & & \multicolumn{2}{c|}{ARMA} & \multicolumn{2}{c|}{MA}\\
 \hline
& & \multicolumn{2}{c|}{$\varrho(\hat{G}||G)$} & \multicolumn{2}{c|}{$\varrho(G||\hat{G})$}\\
\hline
\multirow{2}{*}{$m=N^{1/3}$} & $\hat{f}_0$ & 35.29 (0.020) & 34.74 (0.020) & 13.09 (0.007) & 13.09 (0.007)\\
\cline{2-6}
 & $\frac{1}{2\pi}$ & 32.74 (0.018) & 32.07 (0.018) & 13.11 (0.007) & 13.11 (0.007)\\
 \hline
\multirow{2}{*}{$m=N^{1/4}$} & $\hat{f}_0$ & 31.12 (0.017) & 31.01 (0.017) & 12.84 (0.007) & 12.84 (0.007)\\
\cline{2-6}
 & $\frac{1}{2\pi}$ & 26.39 (0.015) & 26.38 (0.015) & 13.71 (0.008) & 13.71 (0.008)\\
\hline
\end{tabular}
\end{table}
\begin{table}[H]
\small
\centering
\caption{Performance of BIC Criterion of NSCD for ARMA and Invertible MA Processes with Comparison to AutoPARM, WBS, BS, NMCD, and MuBred \label{ARMA2}}
\begin{tabular}{|c|c|c|c|c|c|c|}
\hline
 & \multicolumn{3}{c|}{ARMA} & \multicolumn{3}{c|}{MA}\\
\hline
& $\hat{K}$ & $\varrho(\hat{G}||G)$ & $\varrho(G||\hat{G})$ & $\hat{K}$ & $\varrho(\hat{G}||G)$ & $\varrho(G||\hat{G})$\\
\hline
AutoPARM & 95.3\% &  48.37 (0.027) & 27.62 (0.015) & 99\% & 11.42(0.006) & 11.42 (0.006)\\
\hline
NSCD & 99\% & 31.42 (0.017) & 33.26 (0.018) & 99.6\% & 14.59 (0.008) & 15.69 (0.009)\\
\hline
MuBred & 93.3\% & 40.436 (0.022) & 44.467 (0.025) & 97.7\% & 49.518 (0.028) & 41.444 (0.023) \\
\hline
WBS & 74.9\% & 184.393 (0.102) & 66.288 (0.037) & 91.8\% & 33.586 (0.019) & 45.362 (0.025)\\
\hline
BS & 77.7\% & 188.955 (0.105) & 69.194 (0.038) & 96.3\% & 29.673 (0.016) & 36.611 (0.020)\\
\hline
NMCD & 0\% & 188.64 (0.105) & 324.31 (0.180) & 0.1\% & 160.4434 (0.089) & 326.4282 (0.181)\\
\hline
\end{tabular}
\end{table}
We expect that AutoPARM still works since it is well known that ARMA and invertible MA processes can be approximated by causal AR processes. We see that under ARMA setting, the performance of NSCD is comparable to AutoPARM, while AutoPARM performs better under MA settings. Both NSCD and AutoPARM perform better than WBS and BS. MuBred gives a slightly worse performance since periodogram is not consistent. NMCD fails in both cases as it does in the previous section. The number of change points is overestimated by NMCD again, which is 3.005 and 3.00, respectively.  \par
\subsection{Non-invertible MA Process \label{ch3simu3}}
We will show simulation results when samples are generated from non-invertible MA process. In theory causal AR cannot approximate non-invertible MA process. Results are given in Table \ref{MA1} and \ref{MA2}.\par
Non-invertible MA (Case 4):
\begin{itemize}
\item $X_t=(1+2B+B^2+5B^3)\xi_t$, $1\leq t\leq 500$,
\item $X_t=(1-2B+2B^2-5B^3)\xi_t$, $501\leq t\leq 1100$,
\item $X_t=(1+2B-B^2+5B^3)\xi_t$, $1101\leq t\leq 1800$.
\end{itemize}
\begin{table}[H]
\centering
\caption{Performance of NSCD for Non-Invertible MA Processes with Different Baseline Functions and Bandwidths When $K$ Is Known\label{MA1}}
\begin{tabular}{|c|c|c|c|}
\hline
 &  & \multicolumn{2}{c|}{MA}\\
 \hline
& & $\varrho(\hat{G}||G)$ & $\varrho(G||\hat{G})$\\
\hline
\multirow{2}{*}{$m=N^{1/3}$} & $\hat{f}_0$  & 33.60 (0.019) & 33.60 (0.019)\\
\cline{2-4}
& $\frac{1}{2\pi}$ & 36.00 (0.020) & 36.00 (0.020)\\
 \hline
\multirow{2}{*}{$m=N^{1/4}$} & $\hat{f}_0$ & 40.11 (0.022) & 40.11 (0.022)\\
\cline{2-4}
& $\frac{1}{2\pi}$ & 41.00 (0.023) & 41.00 (0.023)\\
\hline
\end{tabular}
\end{table}

\begin{table}[H]
\centering
\caption{Performance of BIC Criterion of NSCD for Non-Invertible MA Processes with Comparison to AutoPARM, WBS, BS, NMCD, and MuBred \label{MA2}}
\begin{tabular}{|c|c|c|c|c|c|c|}
\hline
 & \multicolumn{3}{c|}{MA}\\
\hline
& $\hat{K}$ & $\varrho(\hat{G}||G)$ & $\varrho(G||\hat{G})$\\
\hline
AutoPARM & 36.7\% & 408.60 (0.227) & 27.24 (0.015)\\
\hline
NSCD & 99.7\% & 37.49 (0.021) & 38.02 (0.021)\\
\hline
MuBred & 8.4\% & 216.203 (0.120) &  26.724 (0.015)\\
\hline
WBS & 18.2\% & 488.674 (0.271) & 90.493 (0.050)\\
\hline
BS & 8.6\% & 576.462 (0.320) & 56.975 (0.032)\\
\hline
NMCD & 0\% & 165.54 (0.092) & 326.701 (0.182)\\
\hline
\end{tabular}
\end{table}
From Table \ref{MA2}, we can see that NSCD works for Case 4, while AutoPARM and MuBred fail since the number of change points are all underestimated. The possible reason is that the marginal variance of Case 4 does not vary much. \par
\subsection{Random Noise without the Existence of Higher Moments}
Next, we investigate the results when $\xi$ does not have higher-order expectations. Here $\xi\sim \frac{1}{\sqrt{2}}t(4)$ so that the variance of $\xi$ is still 1. The results for four cases are shown in Table \ref{highorder}, when $K$ is known.
\begin{table}[H]
\centering
\caption{Performance of NSCD for Random Noise with t(4) \label{highorder}}
\begin{tabular}{|c|c|c|}
\hline
 & $\varrho(\hat{G}||G)$ & $\varrho(G||\hat{G})$\\
\hline
Case 1  & 24.24 (0.012) & 24.37 (0.012) \\
\hline
Case 2 & 29.67 (0.017) & 28.91 (0.016)\\
\hline
Case 3 & 11.74 (0.007) & 11.74 (0.007)\\
\hline
Case 4 & 36.36 (0.020) & 36.36 (0.020)\\
\hline
\end{tabular}
\end{table}

Apparently, from Table \ref{highorder}, it is clear that Assumption 1 can be slightly violated in application, while still achieving good performance.\par
\subsection{Investigation of Smaller Sample Size}
In previous sections, we discuss the performances when sample size is about 2000, which is large. Here we shrink the sample size by half, and investigate the estimate accuracy as well as the choice of bandwidth. We set $ml=200$, and $K_{\mathrm{max}}=5$. Baseline function is chosen to be $\hat{f}_0$. Results are shown in Table \ref{smallsample}.
\begin{table}[H]
\centering
\caption{Performance of NSCD with Different Bandwidths While the Sample Size Is Reduced by Half \label{smallsample}}
\begin{tabular}{|c|c|c|c|}
\hline
 & Bandwidth & $\varrho(\hat{G}||G)$ & $\varrho(G||\hat{G})$\\
\hline
\multirow{2}{*}{Case 1} & $m=N^{1/3}$ & 23.195 (0.023) & 23.263 (0.023)\\
\cline{2-4}
 & $m=N^{1/4}$ & 27.716 (0.027) & 28.075 (0.027) \\
\hline
\multirow{2}{*}{Case 2} & $m=N^{1/3}$ & 31.72 (0.035) & 30.718 (0.035)\\
\cline{2-4}
 & $m=N^{1/4}$ & 28.104 (0.031) & 27.734 (0.030)\\
\hline
\multirow{2}{*}{Case 3} & $m=N^{1/3}$ & 14.664 (0.016) & 14.158 (0.016)\\
\cline{2-4}
 & $m=N^{1/4}$ & 12.652 (0.014) & 12.652 (0.014)\\
\hline
\multirow{2}{*}{Case 4} & $m=N^{1/3}$ & 37.28 (0.041) & 37.272 (0.041)\\
\cline{2-4}
 & $m=N^{1/4}$ & 44.3 (0.05) & 43.877 (0.049)\\
\hline
\end{tabular}
\end{table}
\begin{table}[H]
\centering
\caption{Performance of BIC Criterion of NSCD with Different Bandwidths While The Sample Size Is Reduced by Half \label{smallsamplebic}}
\begin{tabular}{|c|c|c|c|c|}
\hline
 & bandwidth & $\hat{K}$ & $\varrho(\hat{G}||G)$ & $\varrho(G||\hat{G})$\\
\hline
\multirow{2}{*}{Case 1}& $m=N^{1/3}$ & 96\% & 24.972 (0.024) & 30.44 (0.030)\\
\cline{2-5}
 & $m=N^{1/4}$ & 90.4\% & 35.268 (0.040) & 39.494 (0.044)\\
\hline
\multirow{2}{*}{Case 2} & $m=N^{1/3}$ & 96.2\% & 38.326 (0.043) & 28.1 (0.031)\\
\cline{2-5}
 & $m=N^{1/4}$ & 98.7\% & 27.494 (0.031) & 26.384 (0.029)\\
\hline
\multirow{2}{*}{Case 3} & $m=N^{1/3}$ & 96.5\% & 31.458 (0.035) &  12.73 (0.014)\\
\cline{2-5}
 & $m=N^{1/4}$ & 96\% & 32.484 (0.036) & 11.436 (0.013)\\
\hline
\multirow{2}{*}{Case 4} & $m=N^{1/3}$ & 95\% & 41.228 (0.046) & 40.246 (0.045)\\
\cline{2-5}
 & $m=N^{1/4}$ & 92\% & 47.67 (0.053) & 53.254 (0.059)\\
\hline
\end{tabular}
\end{table}
As shown in the table, the estimate accuracy is not affected by the choice of bandwidth when sample size is relatively small. The estimation accuracy of $\kappa$ is worse than the results in Section \ref{ch3simu1} to \ref{ch3simu3}, since the sample size is decreased by half.
\subsection{Further Investigation of BIC Criterion}
\begin{figure}[H]
\centering
\caption{Influence of $m$ and $ml$ on BIC Criterion with $m=N^{1/4}$ and $ml=350$}
\subfloat[AR Processes]{\label{choicec1}\includegraphics[scale=0.25]{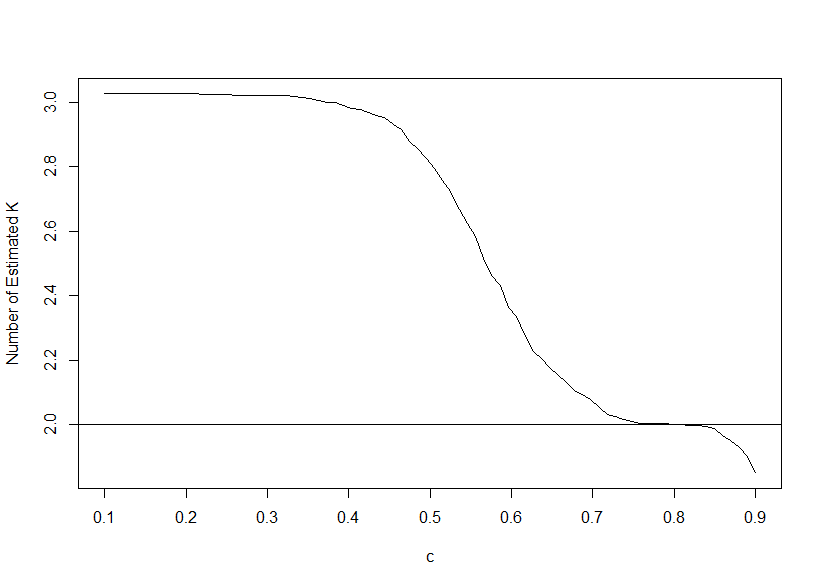}}
\subfloat[ARMA Processes]{\label{choicec2}\includegraphics[scale=0.25]{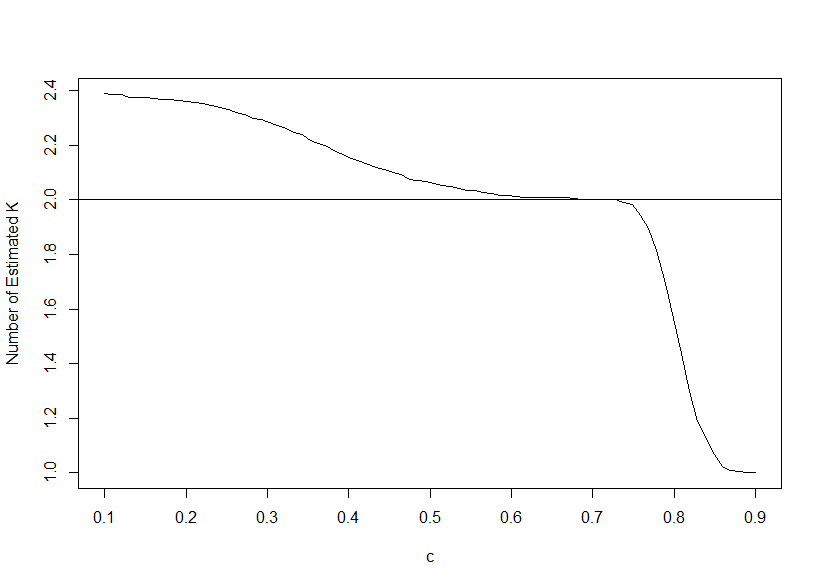}}\\
\subfloat[Invertible MA Processes]{\label{choicec3}\includegraphics[scale=0.25]{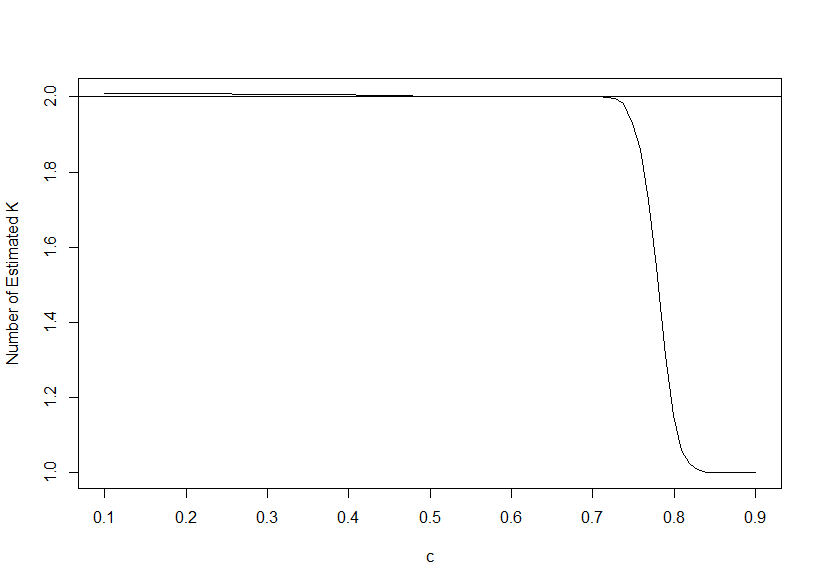}}
\subfloat[Non-Invertible MA Processes]{\label{choicec4}\includegraphics[scale=0.25]{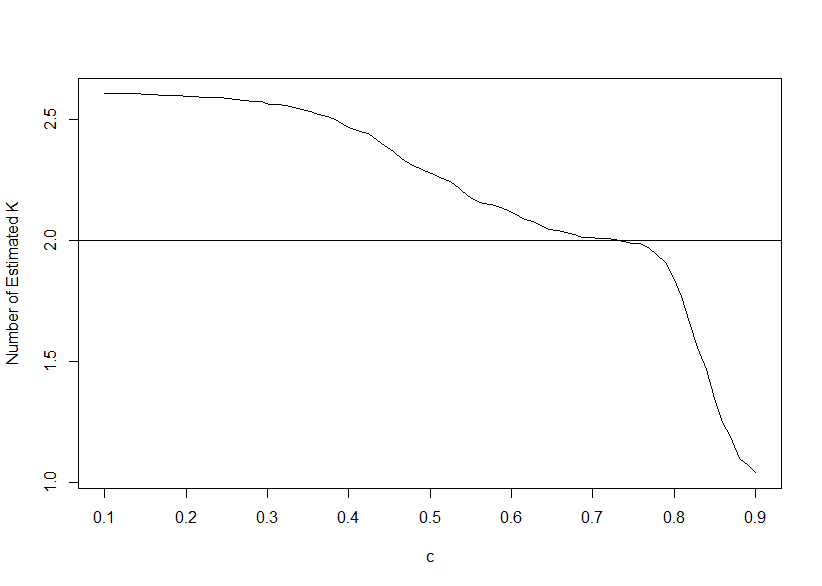}}
\end{figure}

\begin{figure}[H]
\centering
\caption{Influence of $m$ and $ml$ on BIC Criterion with $m=N^{1/3}$ and $ml=350$}
\subfloat[AR Processes]{\label{choicec5}\includegraphics[scale=0.25]{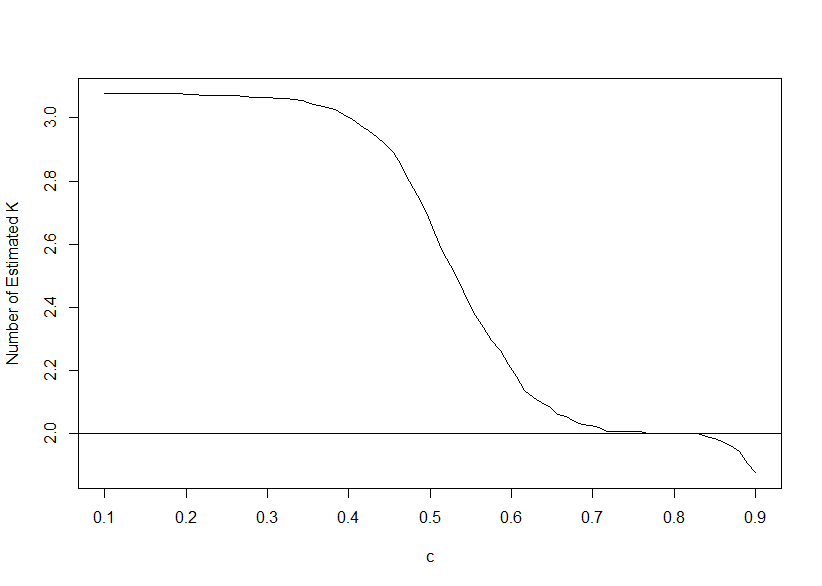}}
\subfloat[ARMA Processes]{\label{choicec6}\includegraphics[scale=0.25]{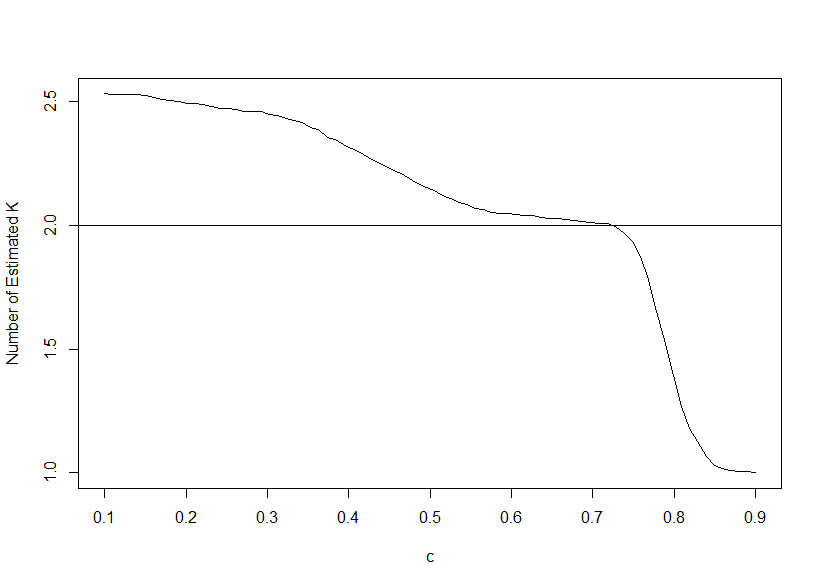}}\\
\subfloat[Invertible MA Processes]{\label{choicec7}\includegraphics[scale=0.25]{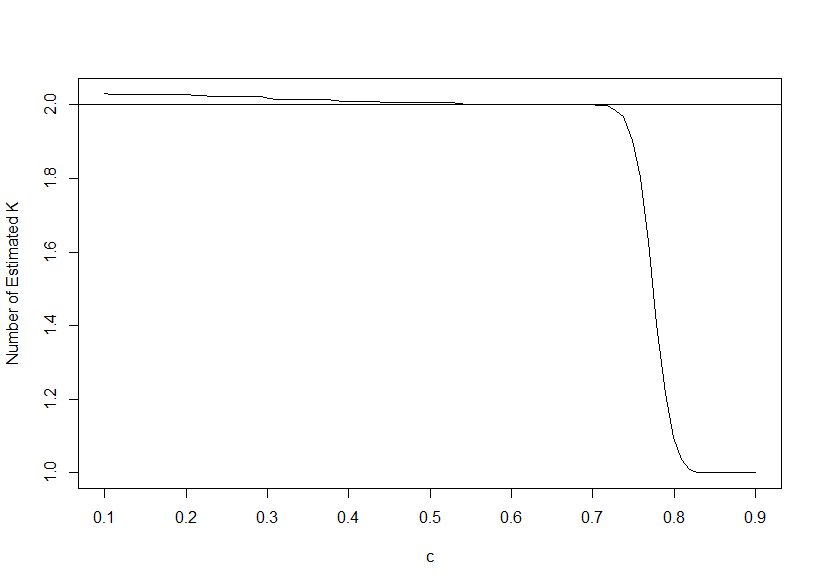}}
\subfloat[Non-Invertible MA Processes]{\label{choicec8}\includegraphics[scale=0.25]{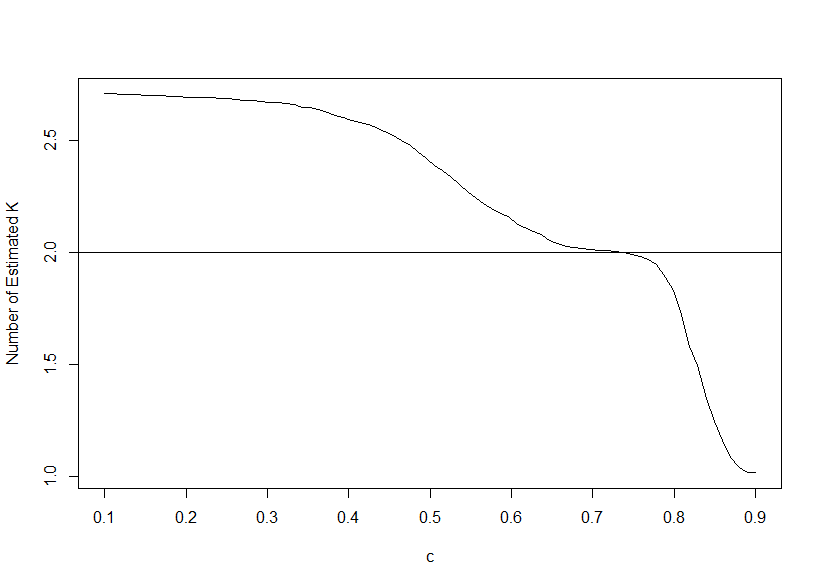}}
\end{figure}

In this section, we investigate the performance with different choice of $N^c$. $c$ ranges from 0.1 to 0.9 and $\hat{K}$ will be plotted in all the cases mentioned above with two choices of bandwidth. All the results are shown from Figures \ref{choicec1} to \ref{choicec8} with baseline function $\hat{f}_0$. From figures, we can see that $c=0.73$ is a good choice for all the four cases.\par
\begin{figure}[H]
\centering
\caption{Influence of $m$ and $ml$ on BIC Criterion with $m=N^{1/3}$ and $ml=300$}
\subfloat[AR Processes]{\label{investbic300p3c1}\includegraphics[scale=0.25]{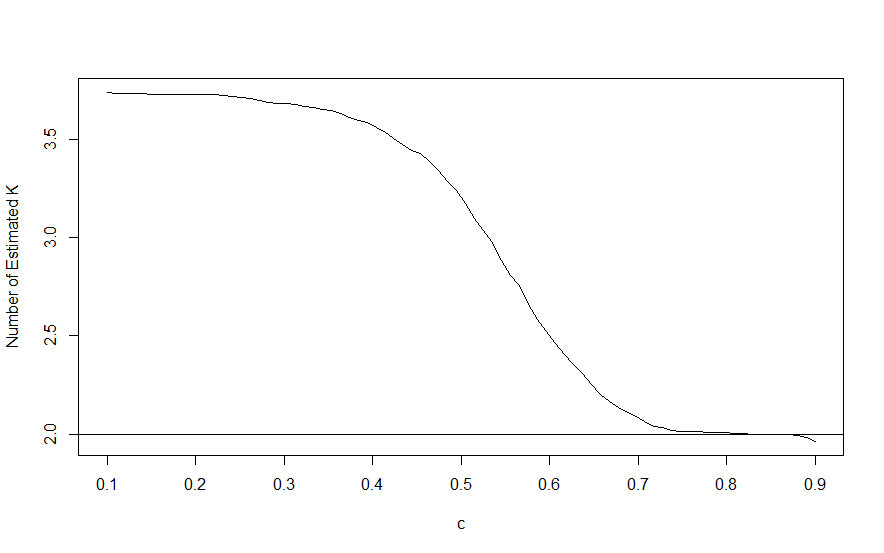}}
\subfloat[ARMA Processes]{\label{investbic300p3c2}\includegraphics[scale=0.25]{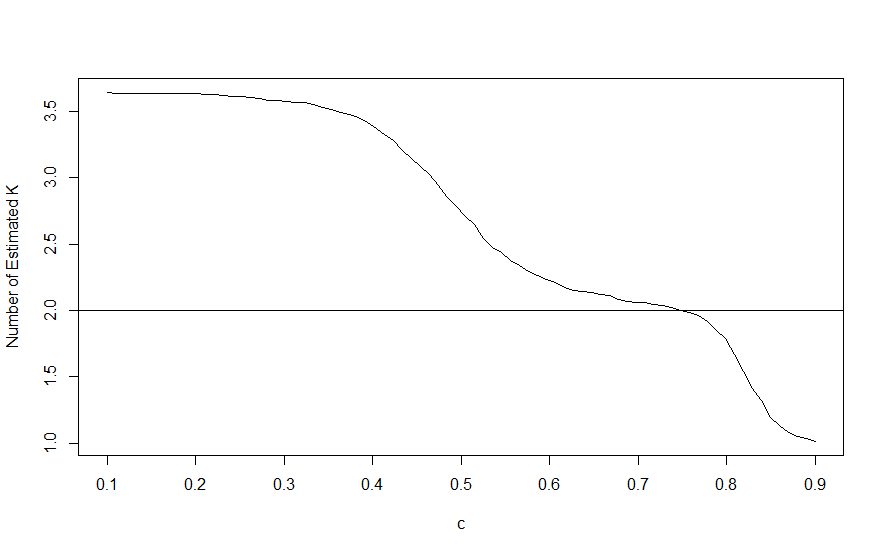}}\\
\subfloat[Invertible MA Processes]{\label{investbic300p3c3}\includegraphics[scale=0.25]{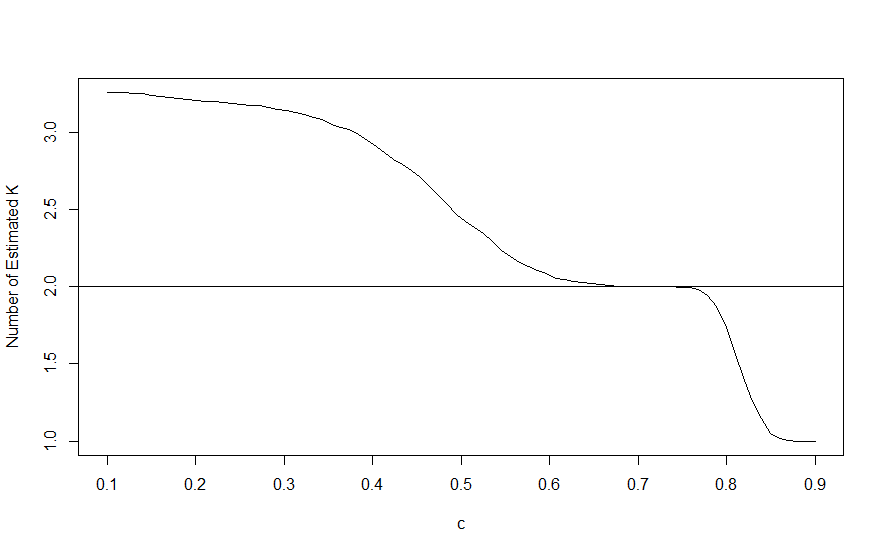}}
\subfloat[Non-Invertible MA Processes]{\label{investbic300p3c4}\includegraphics[scale=0.25]{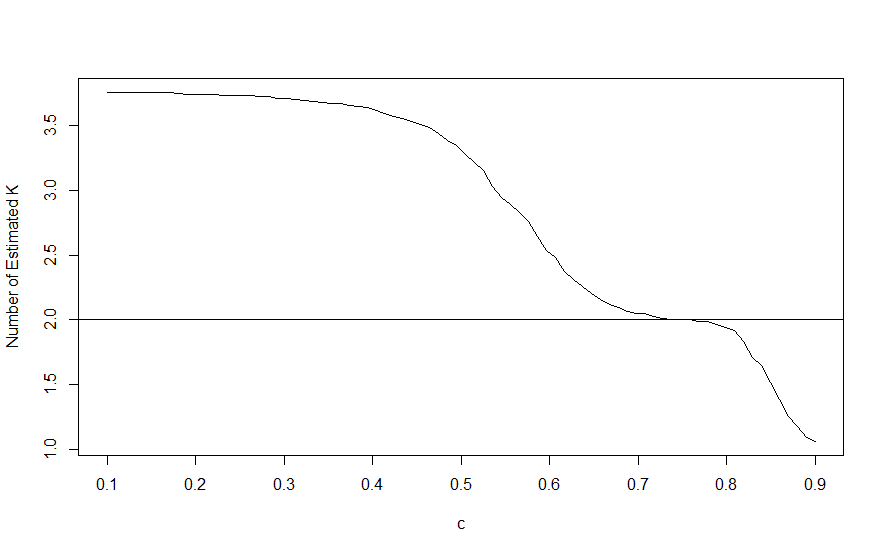}}
\end{figure}
\begin{figure}[H]
\centering
\caption{Influence of $m$ and $ml$ on BIC Criterion with $m=N^{1/4}$ and $ml=300$}
\subfloat[AR Processes]{\label{investbic300p4c1}\includegraphics[scale=0.25]{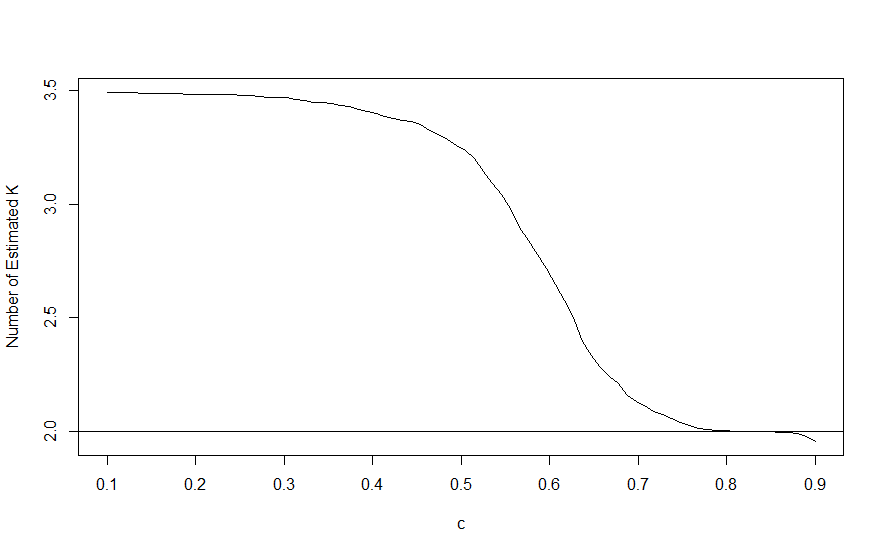}}
\subfloat[ARMA Processes]{\label{investbic300p4c2}\includegraphics[scale=0.25]{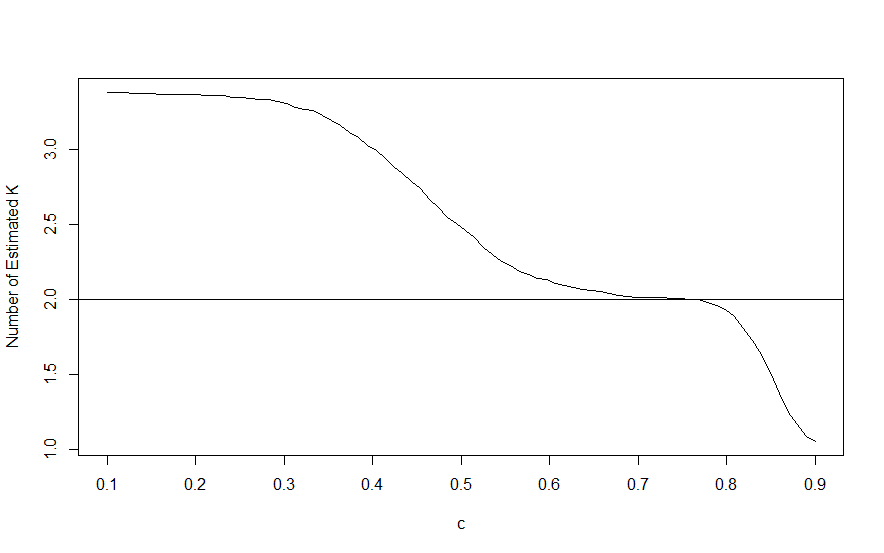}}\\
\subfloat[Invertible MA Processes]{\label{investbic300p4c3}\includegraphics[scale=0.25]{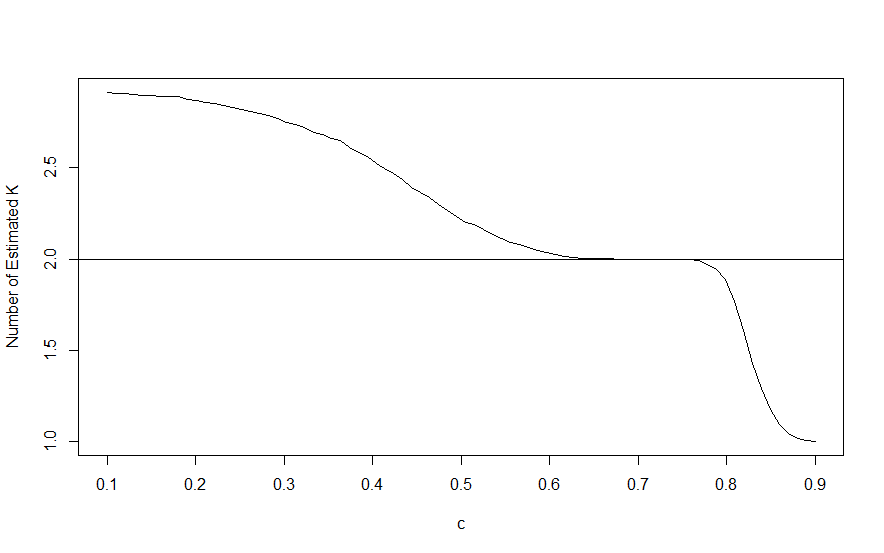}}
\subfloat[Non-Invertible MA Processes]{\label{investbic300p4c4}\includegraphics[scale=0.25]{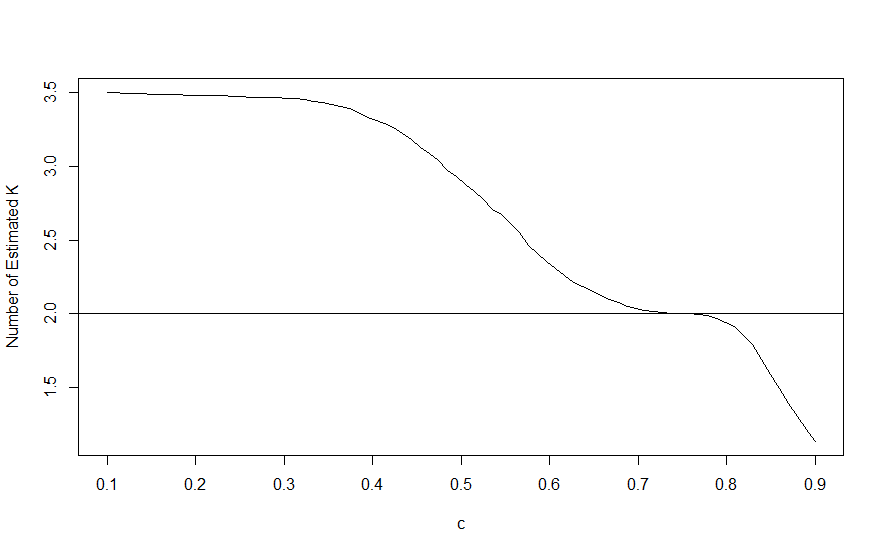}}
\end{figure}
\begin{figure}[H]
\centering
\caption{Influence of $m$ and $ml$ on BIC Criterion with $m=N^{1/3}$ and $ml=250$}
\subfloat[AR Processes]{\label{investbic250p3c1}\includegraphics[scale=0.25]{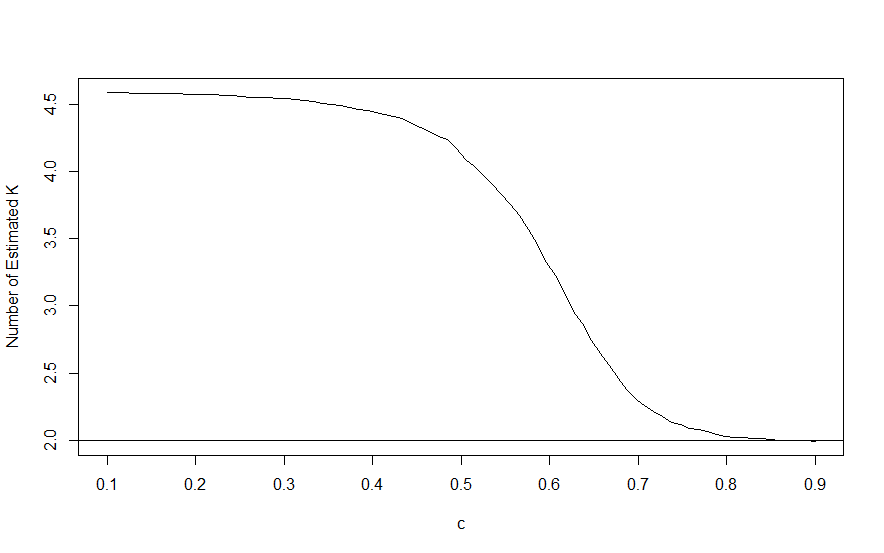}}
\subfloat[ARMA Processes]{\label{investbic250p3c2}\includegraphics[scale=0.25]{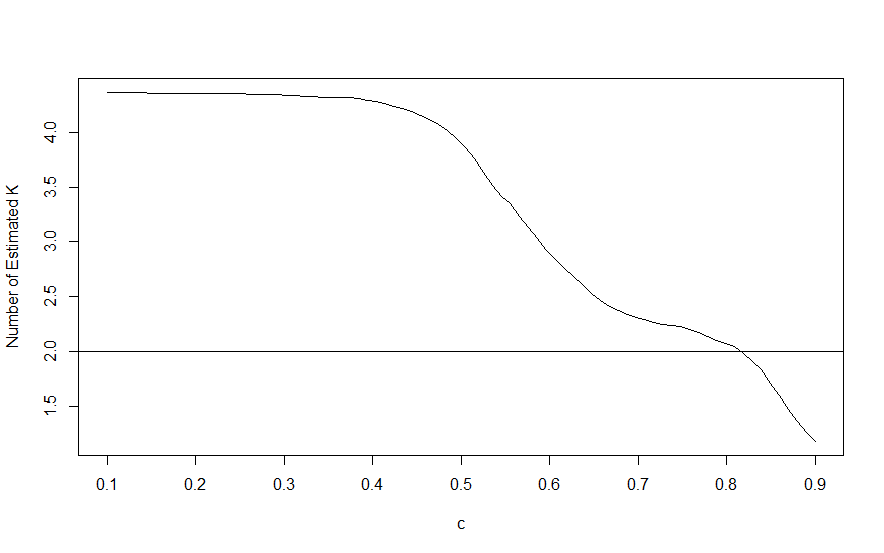}}\\
\subfloat[Invertible MA Processes]{\label{investbic250p3c3}\includegraphics[scale=0.25]{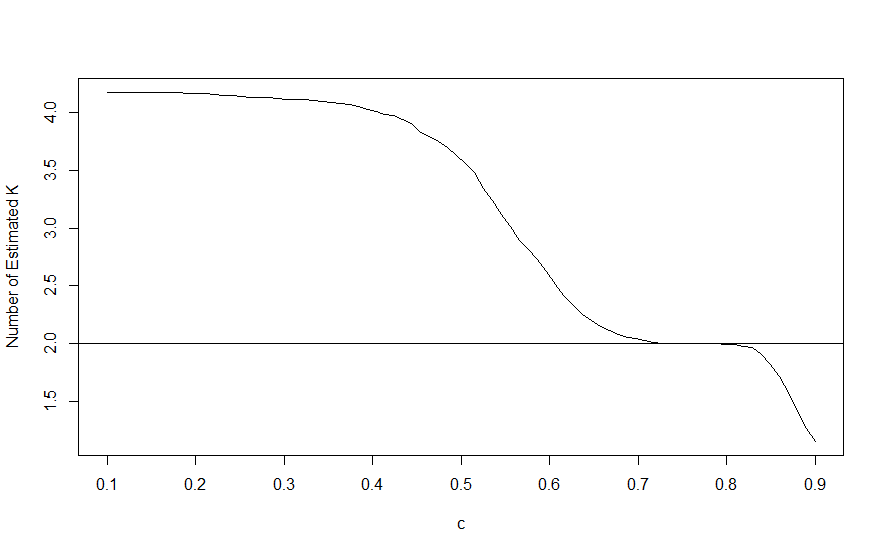}}
\subfloat[Non-Invertible MA Processes]{\label{investbic250p3c4}\includegraphics[scale=0.25]{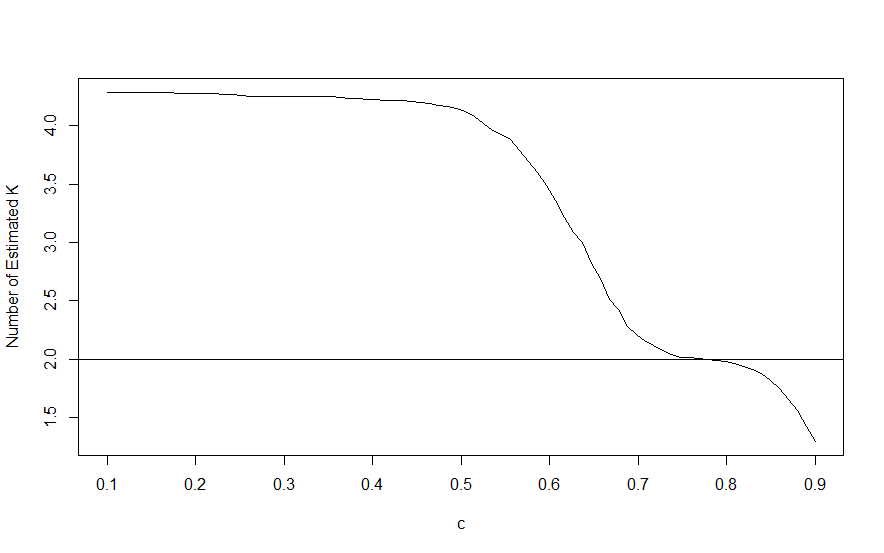}}
\end{figure}
\begin{figure}[H]
\centering
\caption{Influence of $m$ and $ml$ on BIC Criterion with $m=N^{1/4}$ and $ml=250$}
\subfloat[AR Processes]{\label{investbic250p4c1}\includegraphics[scale=0.25]{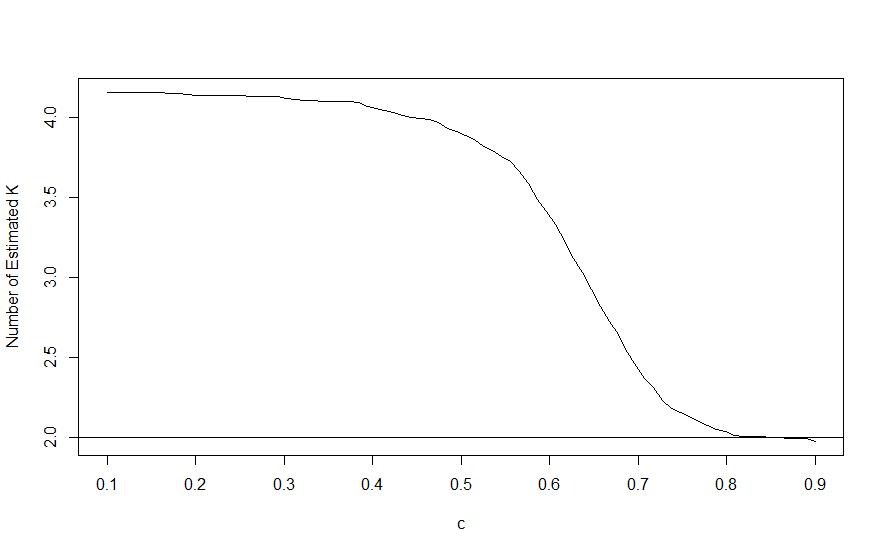}}
\subfloat[ARMA Processes]{\label{investbic250p4c2}\includegraphics[scale=0.25]{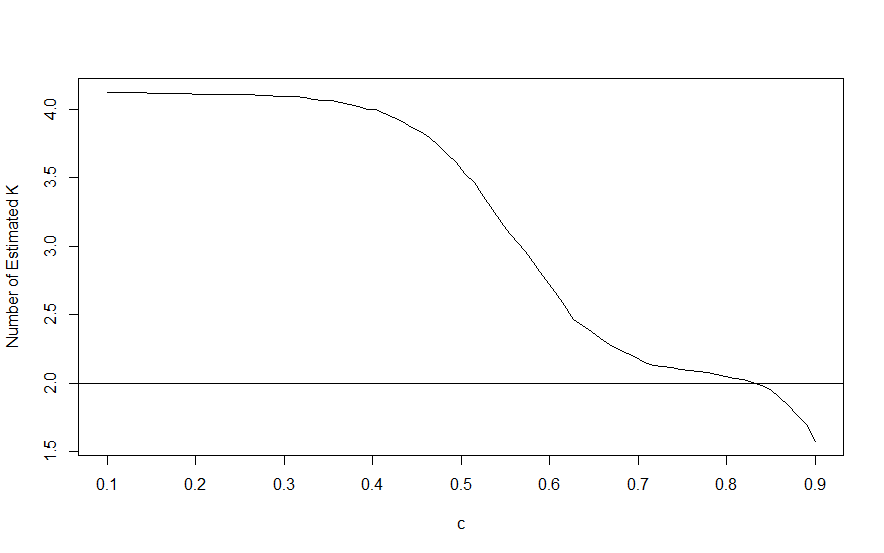}}\\
\subfloat[Invertible MA Processes]{\label{investbic250p4c3}\includegraphics[scale=0.25]{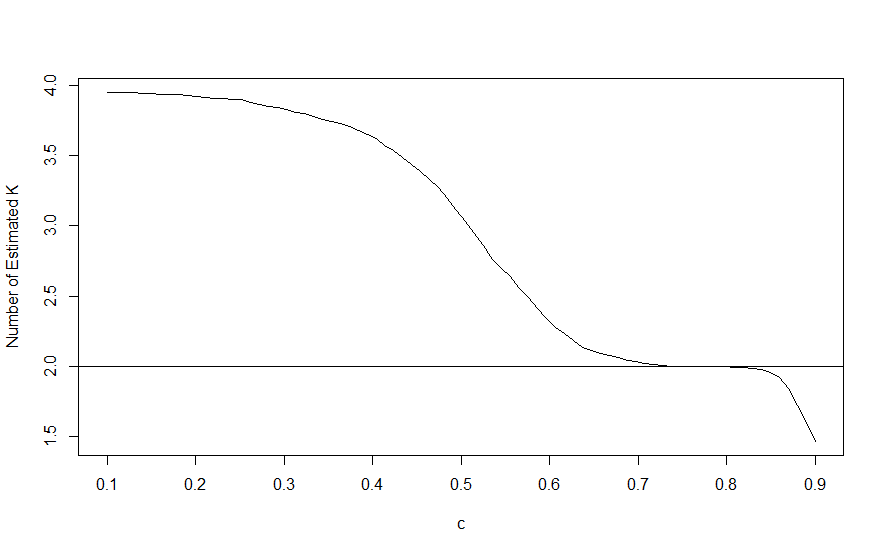}}
\subfloat[Non-Invertible MA Processes]{\label{investbic250p4c4}\includegraphics[scale=0.25]{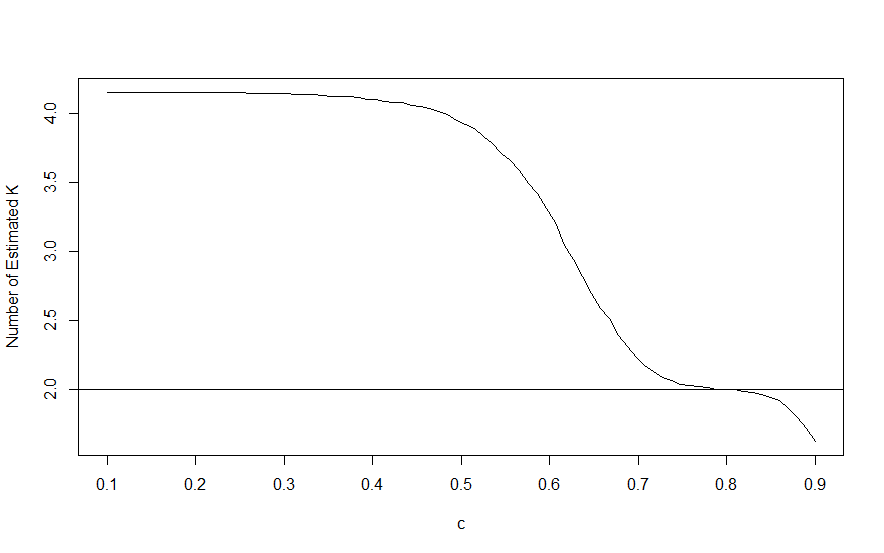}}
\end{figure}

Next, we are going to simulate the performances of BIC when $ml$ changes. Since from the previous section,we can see that the choice of $C_N$ depends on $ml$. Figures \ref{investbic300p3c1} to \ref{investbic300p4c4} show the results when $ml=300$, while Figures \ref{investbic250p3c1} to \ref{investbic250p4c4} give the performances when $ml=250$.\par
As we can see from the figures, $c=0.73$ is not a good choice in general, and it will overestimate the number of change points. Although overestimating is tolerable since we do not want to miss the true change points, we still suggest that a sufficiently large choice for $ml$ is necessary. Since the maximum possible number of change points is $[N/ml]$. We suggest that one should choose a reasonable $K_{\mathrm{max}}$ depending on some prior information, then choose a $ml$ satisfying $ml\leq N/K_{\mathrm{max}}$.\par

\clearpage

\clearpage

\subsection{Influence of Change Point Searching Unit}
In this section, we investigate the influence of $n_{\mathrm{su}}$. For all the four cases, $n_{\mathrm{su}}=10$, which means that for Case 2-4, $\tau_1^0,\ldots \tau_K^0$ are all a multiple of $n_{\mathrm{su}}$, while for Case 1, the true change points are not divisible by $n_{\mathrm{su}}$. We set $m=N^{1/4}$ and baseline function $f=\hat{f}_0$ for four cases. The results are shown in Table \ref{su1} and \ref{su2}.
\begin{table}[H]
\centering
\caption{Performance of NSCD with Searching Unit When $K$ Is Known \label{su1}}
\begin{tabular}{|c|c|c|c|}
\hline
 & $\varrho(\hat{G}||G)$ & $\varrho(G||\hat{G})$\\
\hline
Case 1 & 22.9 (0.011) & 22.9 (0.011)\\
\hline
Case 2 & 29.26 (0.016) & 29.26 (0.016)\\
\hline
Case 3 & 10.14 (0.006) & 10.14 (0.006)\\
\hline
Case 4 & 40.52 (0.023) & 40.52 (0.023)\\
\hline
\end{tabular}
\end{table}
\begin{table}[H]
\centering
\caption{Performance of BIC Criterion of NSCD with Searching Unit\label{su2}}
\begin{tabular}{|c|c|c|c|}
\hline
 & $\hat{K}$ & $\varrho(\hat{G}||G)$ & $\varrho(G||\hat{G})$\\
\hline
Case 1 &  96.2\% & 22.696 (0.011) & 39.72 (0.019)\\
\hline
Case 2 & 97.2\% & 42.02 (0.0233) & 30.68 (0.017)\\
\hline
Case 3 & 99.6\% & 12.54 (0.007) & 10.14 (0.006)\\
\hline
Case 4 & 98.6\% & 40.56 (0.023) & 43.74 (0.024)\\
\hline
\end{tabular}
\end{table}
We can see that setting an $n_{su}$ will not affect our results much even if the true change points are not a multiple of $n_{su}$ in Case 1, so it is safe to apply this in application, which could boost the computation as well as give accurate estimations.

\section{Case Study \label{ch3s5}}
\subsection{Simulated Data}
\begin{figure}[H]
\centering
\caption{Estimated Change Points by NSCD for Different Cases}
\subfloat[AR Processes]{\label{ARdata}\includegraphics[scale=0.25]{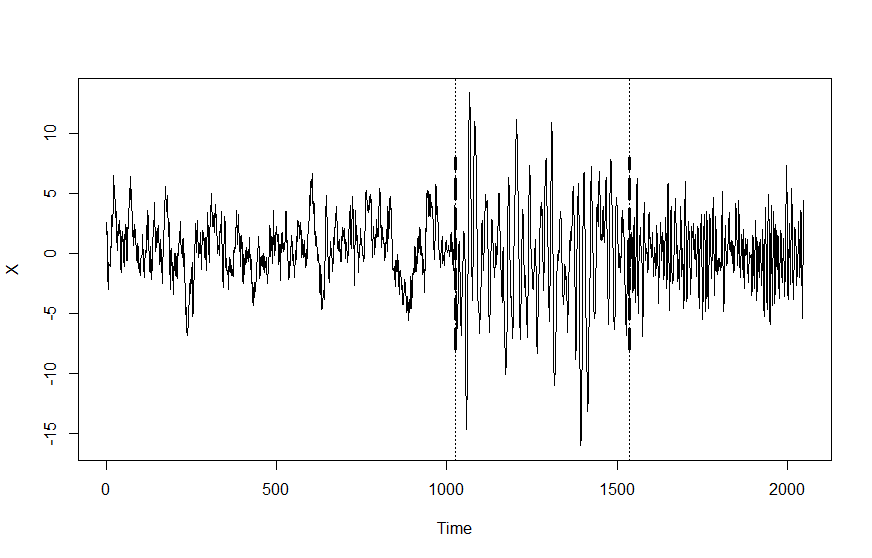}}
\subfloat[ARMA Processes]{\label{ARMAdata}\includegraphics[scale=0.25]{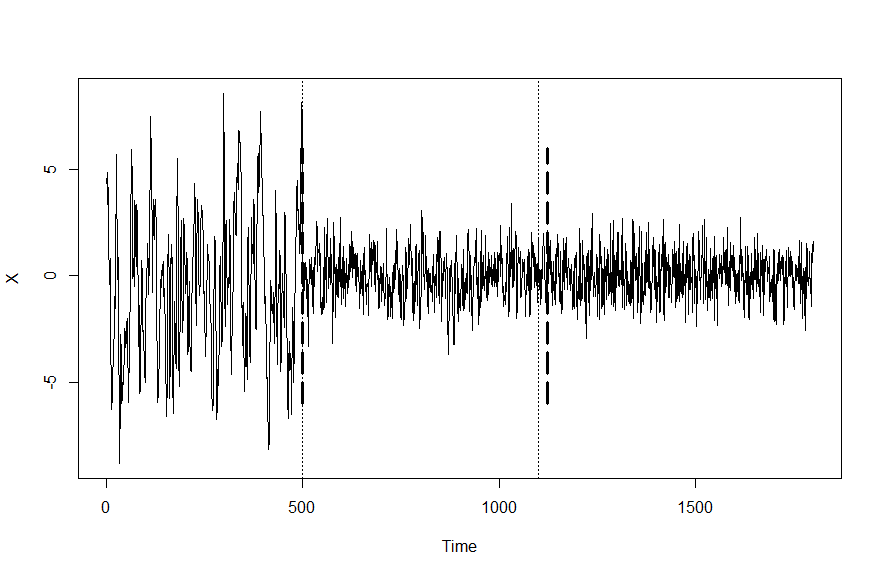}}\\
\subfloat[Invertible MA Processes]{\label{MAdata}\includegraphics[scale=0.22]{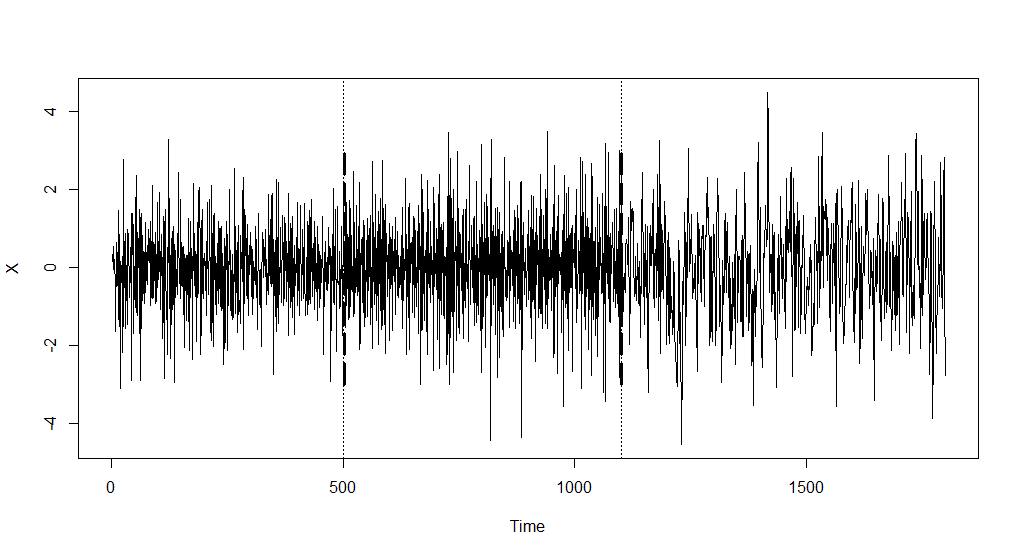}}
\subfloat[Non-Invertible MA Processes]{\label{nonMAdata}\includegraphics[scale=0.26]{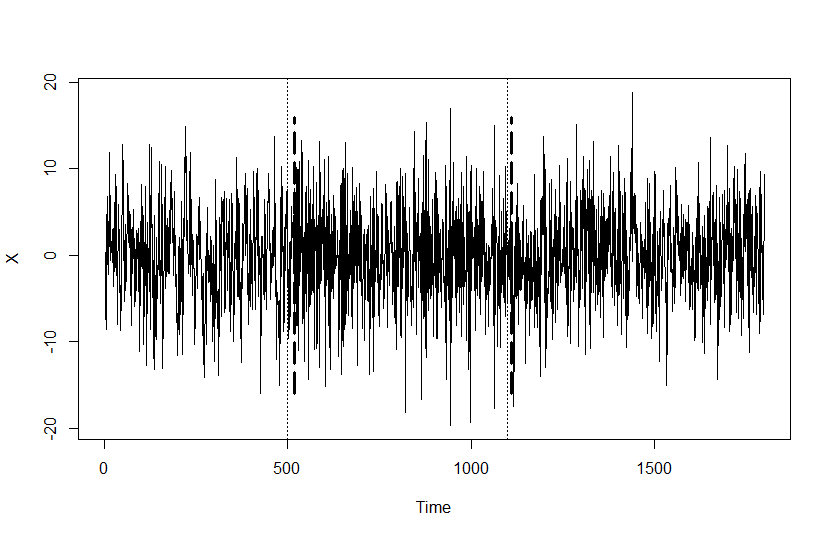}}
\end{figure}

\begin{figure}[H]
\centering
\caption{Estimated Spectrums by NSCD for Different Cases}
\subfloat[AR Processes]{\label{ARspectrum}\includegraphics[scale=0.25]{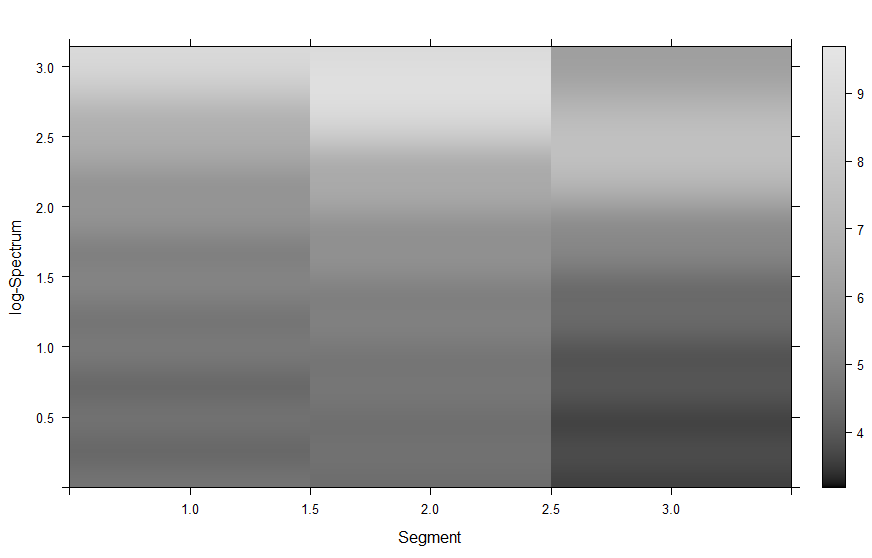}}
\subfloat[ARMA Processes]{\label{ARMAspectrum}\includegraphics[scale=0.25]{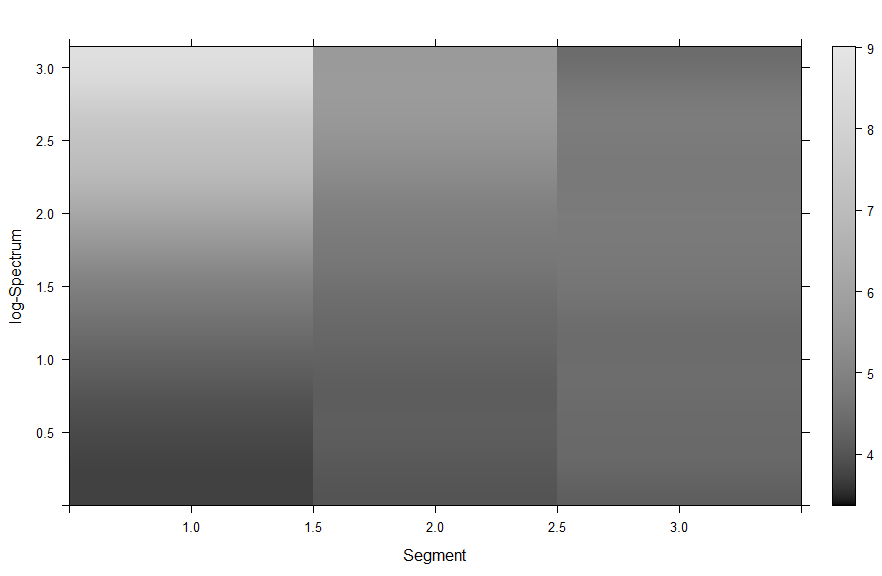}}\\
\subfloat[Invertible MA Processes]{\label{MAspectrum}\includegraphics[scale=0.23]{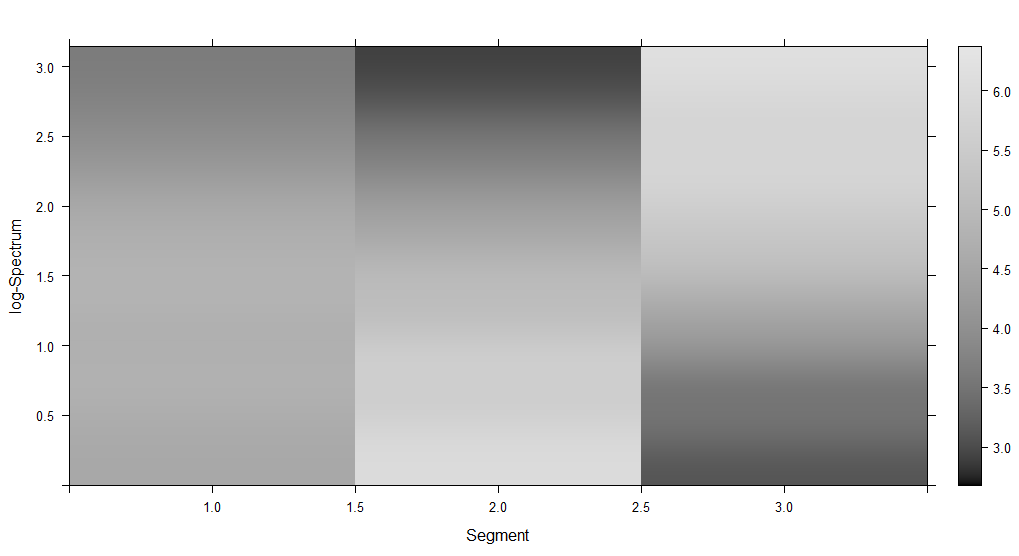}}
\subfloat[Non-Invertible MA Processes]{\label{nonMAspectrum}\includegraphics[scale=0.25]{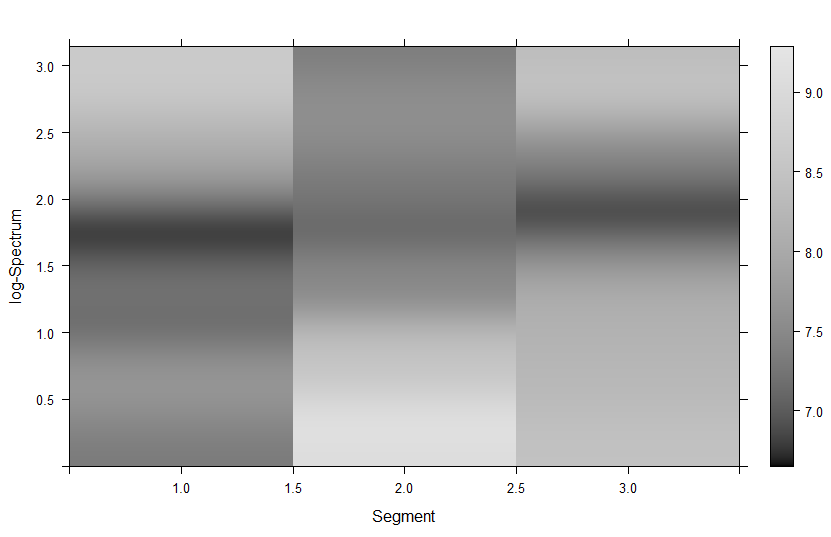}}
\end{figure}

In the simulation section, we check the performances of our method based on 4 cases, and compare NMSD with other methods. Here we will investigate the change point detection of these 4 cases more directly. Figure \ref{ARdata}-\ref{nonMAdata} are the realizations of Case 1 to Case 4. The long-dashed lines are the locations of estimated change points while the dotted lines represent true change points. It is easy to see from Figure \ref{ARdata} that the existence of change points are obvious, since marginal variance of the second segment is bigger. In Case 2, the second change point is far less obvious than the first one which can be seen in Figure \ref{ARMAdata}. In Figure \ref{MAdata}, since the observations in the first two segments concentrate more around their mean compared to the third segment, the second change point may be captured by eyes. In Figure \ref{nonMAdata}, two change points are not apparent any more. However, from the perspective of spectrum, those change points can be easily detected. Figure \ref{ARspectrum}-\ref{nonMAspectrum} give the estimated spectral density functions. In Case 4, we can see that the power of spectrum of the first segment concentrates more at higher and lower frequencies, while the spectrum of the second part has more power at low frequency. The spectrum of the third part is similar to the first one, but it has more power in the low frequency range.

\subsection{Electroencephalography Recordings for Seizure}
\begin{figure}[H]
\centering
\caption{Analysis of Different Channels of EEG Data}
\subfloat[Channel FP1-F3]{\label{EEGF3}\includegraphics[scale=0.24]{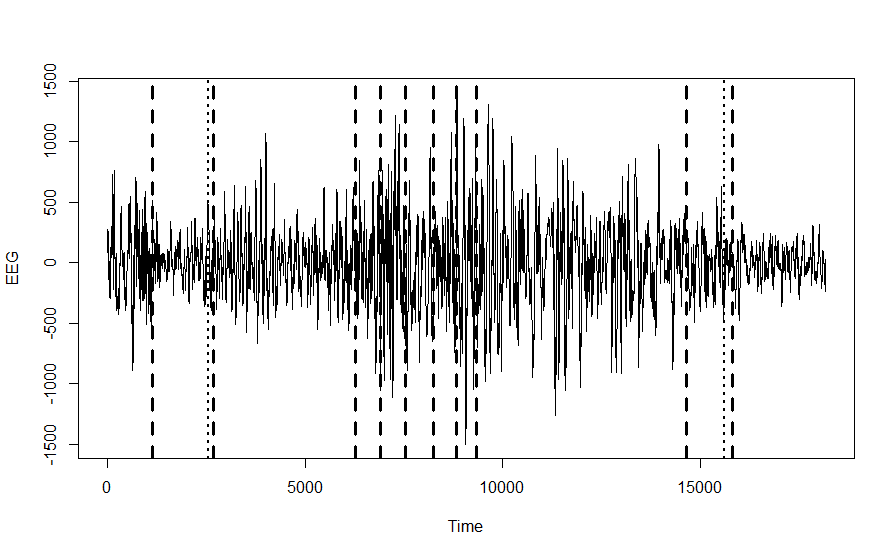}}
\subfloat[Channel FP1-F7]{\label{EEGF7}\includegraphics[scale=0.24]{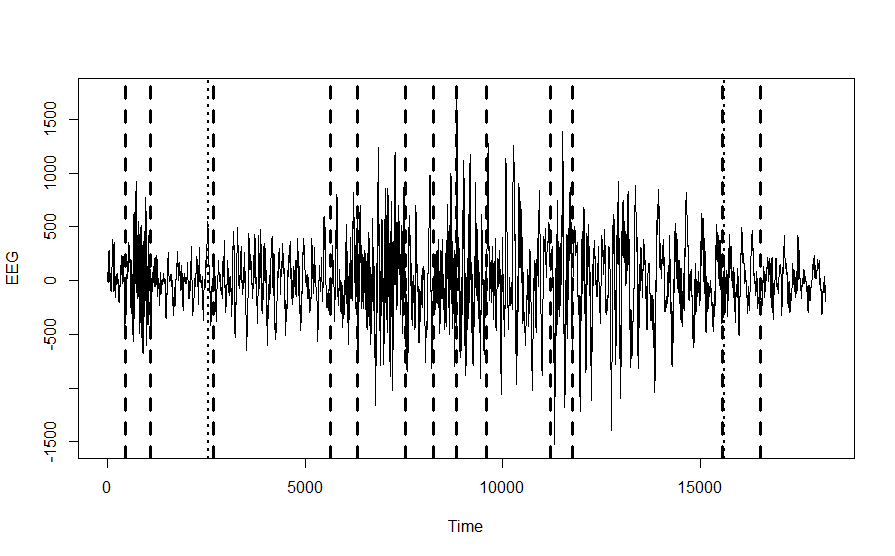}}
\end{figure}

In this section, we investigate the performance of NSCD when applied to EEG data for seizure (Goldberger et al. 2000). The data is retrieved from https://www.physionet.org/pn6/chbmit/. For each subject, the EEG signal was recorded into several data files, which was one-hour long. Here all subjects were monitored for several days to trace their states, so we only analyze the data file with seizures. The sampling rate is 256 Hz, so the number of observations in each recording is 9216000, which is huge. In the data file we choose (which is chb01\_16), seizure happened only once with duration 51 seconds, so the number of change points is 2. The duration of seizure is very short compared with the total length of this recording, so we analyze a particular part of the data which begins at 10 seconds before the seizure and ends at the 10 seconds after the seizure, so the sample size is 18176, which is still large. So here we will set a change point searching unit equal to 64, which is 0.25 second, to ease the computation burden. What is more, the minimal length $n_{\mathrm{min}}$ is 256 which is 1 second. There are totally 23 channels in EEG recording, we apply our method on channel ``FP1-F3'' and ``FP1-F7''. The results are shown in Figure \ref{EEGF3} and \ref{EEGF7}. The vertical dotted lines represent the locations of the beginning and ending of seizure, while the vertical dashed lines are the estimated change points. To demonstrate the changes in spectrums, we plot the estimated spectrums in Figure {F3s} and {F7s}. As we can see, NSCD can successfully detect the true change points, while almost all the other estimates fall between the true change points. This is because during seizure, the EEG recordings change abruptly, which will bring more non-stationarity into the time series.\par

\begin{figure}[H]
\centering
\caption{Estimated Spectrums of Different Channels of EEG Data}
\subfloat[Channel FP1-F3]{\label{F3s}\includegraphics[scale=0.24]{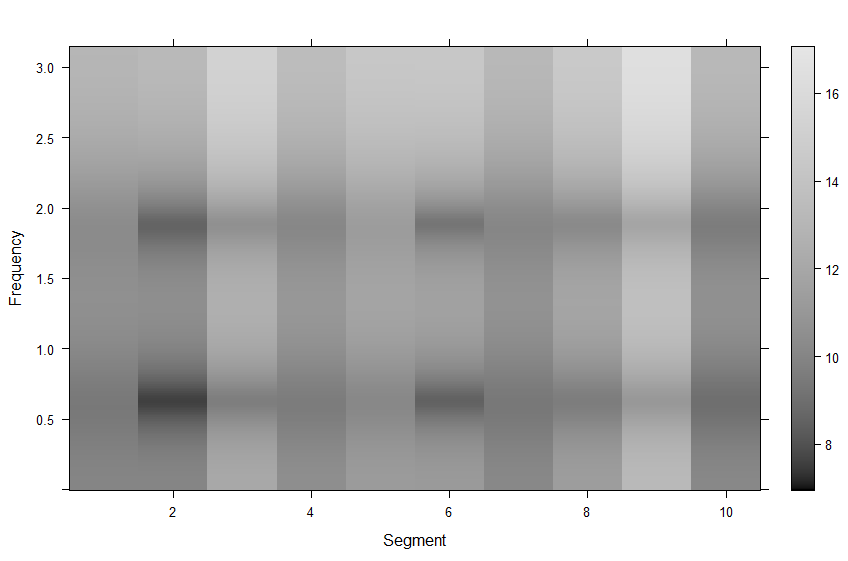}}
\subfloat[Channel FP1-F7]{\label{F7s}\includegraphics[scale=0.24]{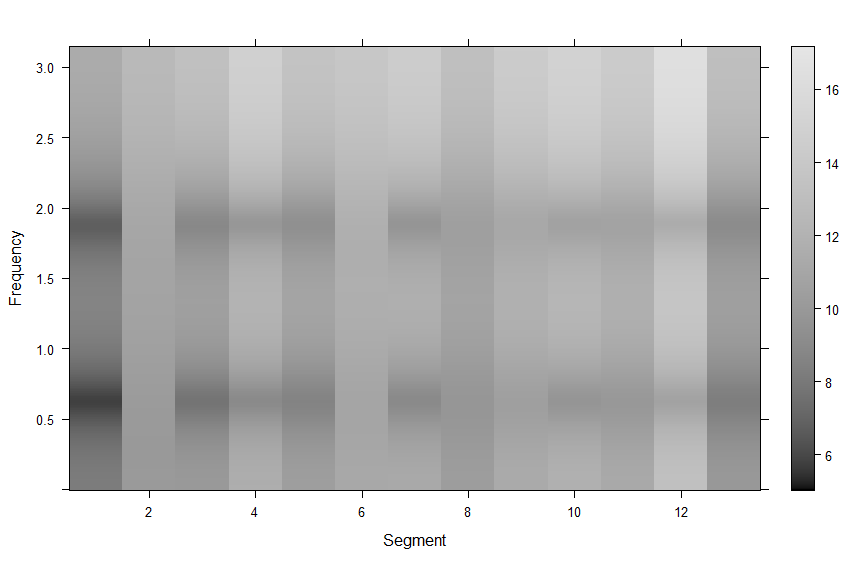}}
\end{figure}

\section{Conclusion}
In this article, we propose a change point detection method based on spectral density functions for non-stationary time series. We assume that non-stationary time series can be segmented into several linear processes. Then Kullback-Leibler divergence is applied to measure the discrepancy between different spectral density functions. A BIC criterion is suggested to estimate the number of change points. Due to the separable structure of objective function, we use Dynamic Programming to find the estimators. We also show the consistency of our estimators in theory and the estimate accuracy by simulations.

\appendix

\makeatletter   
 \renewcommand{\@seccntformat}[1]{APPENDIX~{\csname the#1\endcsname}.\hspace*{1em}}
 \makeatother

\section{Proofs}
In the appendix, $B_1$ to $B_8$ are appropriate constant and $B_{7\epsilon}$, $B_{8\epsilon}$ are constant with respect to $\epsilon$.\par
\begin{lemma}\label{ch3l1} Suppose $X_t=\sum\limits_{j=-\infty}^{+\infty}a_j\xi_{t-j}$, $1\leq t\leq N$, where $\sum\limits_j|a_j|<+\infty$. $\xi_j\overset{\mathrm{i.i.d}}{\sim}(0,\sigma^2)$. Denote $f(\lambda)$ as the spectral density function of $X_t$. $\hat{f}$ is the estimated spectral density function. Then under Assumptions 1-6,
\begin{center}
$\displaystyle \left|\frac{\hat{f}(\lambda)-f(\lambda)}{f(\lambda)}\right|^{2q}=O_p(\frac{m^{2q}}{N^{q}})$\\
$\displaystyle \max\limits_{\lambda_i\in \Lambda}\left|\frac{\hat{f}(\lambda_i)-f(\lambda_i)}{f(\lambda_i)}\right|^{2q}=O_p(\frac{m^{2q}}{N^{q-1}})$\\
\end{center}
\end{lemma}
\noindent\textbf{Proof}: Following the method of Woodroofe and Van Ness (1967), we separate our proofs into 3 parts.\par
\noindent Part 1: we show that for given $\lambda\in[-\pi,\pi]$, $\displaystyle E|g_N(\lambda)-1|^{2q}=O_p\left(\frac{m^{2q}}{N^q}\right)$. Here $g_N(\lambda)$ is the smoothed spectral density estimation for $\{\xi_1,\cdots,\xi_N\}$, $\xi$
\begin{eqnarray*}
\left| g_N(\lambda)-1\right|^{2q}&=&\left(m\int W(m(u-\lambda))\frac{1}{N}\left|\sum\limits_{t=1}^N\xi_te^{-itu}\right|^2du-\sigma^2\right)^{2q}\\
&=&\frac{1}{N^{2q}}\left(\int W(u)\left(\left(\sum\limits_{t=1}^N\xi_t e^{-it(um^{-1}+\lambda )}\right)^2-N\sigma^2\right)du\right)^{2q}\\
&=&\frac{1}{N^{2q}}\left(2\sum\limits_{t=1}^N\sum\limits_{v=1}^{N-t}\xi_{t}\xi_{t+v}\cos(v\lambda)w(vm^{-1})+\sum\limits_{t=1}^N\xi^2_t-N\sigma^2\right)^{2q}\\
&=&\frac{1}{N^{2q}}\left(Z_N(\lambda)+r_N(\lambda)+r_N\right)^{2q},
\end{eqnarray*}
where $Z_N(\lambda)=\sum\limits_{t=1}^NZ_{N,t}(\lambda)$,
\begin{center}
$Z_{N,t}(\lambda)=2\sum\limits_{v=1}^{m-1}\xi_t\xi_{t+s}w(vm^{-1})\cos(v\lambda)$, \\ $r_N(\lambda)=2\sum\limits_{t=N-m+2}^N\sum\limits_{v=N-t+1}^{m-1}\xi_t\xi_{t+v}w(vm^{-1})\cos(v\lambda)$, \\ $r_N=\sum\limits_{t=1}^N(\xi_t^2-\sigma^2)$.
\end{center}
Since for any real number $b_1,\ldots,b_{p}$, $\left(b_1+\ldots+b_{p}\right)^{2q}\leq 2^{2pq}(b_1^{2q}+\ldots+b_p^{2q})$. So we only need to prove that $\left(Z_N(\lambda)\right)^{2q}$, $\left(r_N(\lambda)\right)^{2q}$, $\left(r_N\right)^{2q}$ are $O_p(N^qm^{2q})$. By Markov's inequality, it suffices to show that $E\left(Z_N(\lambda)\right)^{2q}$, $E\left(r_N(\lambda)\right)^{2q}$, $E\left(r_N\right)^{2q}$ are $O(N^qm^{2q})$.\par
After expanding $\left(r_N(\lambda)\right)^{2q}$, we have
\begin{eqnarray*}
\left(r_N(\lambda)\right)^{2q}&=&\sum\limits_{p=1}^{2q}\sum\limits_{t_1+\ldots+t_p=2q}\sum\limits_{j_1\cdots j_p}\sum\limits_{v_1\cdots v_p}\left(\xi_{j_1}\xi_{j_1+v_1}\right)^{t_1}\cdots\left(\xi_{j_p}\xi_{j_p+v_p}\right)^{t_p}\\
& &w^{t_1}(v_1m^{-1})\cdots w^{t_p}(v_pm^{-1})\cos^{t_1}(v_1\lambda)\cdots\cos^{t_p}(v_p\lambda),
\end{eqnarray*}
where $\{j_1,\ldots,j_p\}\subset\{N-m+2,\ldots,N\}$. If $t_k\geq 2$ for any $k$, then expectation of $\left(\xi_{j_1}\xi_{j_1+v_1}\right)^{t_1}\cdots\left(\xi_{j_p}\xi_{j_p+v_p}\right)^{t_p}$ is not zero. If one of $t_k$ is 1, for example, $t_{k_0}=1$, since $j_k\leq N < j_k+v_k$ for any $k$, it is impossible to find $\xi_{j_{k_0}}$ in $\{\xi_{j_1+v_1},\ldots,\xi_{j_p+v_p}\}$, so the expectation would be zero. To sum up, for all non-zero terms in $E\left(r_N(\lambda)\right)^{2q}$, $t_k$ should be greater than 2. Since for any $s\leq 2q$, $\sum\limits_{k=1}^s t_k=2q$, $w^j(u)\cos^j(u)$, $E\xi_1^{t_1}\cdots\xi_s^{t_s}$ can be bounded by a real number, denoted by $B_1$, so we only need to count the number of non-zero terms. When $t_k=2$ for any $k$, then $p=q$, and the number of non-zero terms are no more than $\left(C_{2q}^2C_{2q-2}^2\cdots C_2^2 C_{m}^q\right)^2=O(m^{2q})$. When at least one of $t_k>2$, say $t_{k_0}=3$, the number of non-zero terms is no more than $\left(C_{2q}^2C_{2q-2}^2\cdots C_4^3 C_{m}^{q-1}\right)^2=O(m^{2q-2})$. So
\begin{center}
$E\left(r_N(\lambda)\right)^{2q}\leq B_1m^{2q}$.
\end{center}
For
\begin{center}
$E(r_N)^{2q}=E(\sum\limits_{t=1}^N(\xi_t^2-\sigma^2))^{2q}=\sum\limits_{p=1}^{2q}\sum\limits_{t_1+\ldots+t_p=2q}\sum\limits_{j_1\cdots j_p}(\xi_{j_1}^2-\sigma^2)^{t_1}\ldots(\xi_{j_p}^2-\sigma^2)^{t_p}\leq B_1 m^{2q}$, \end{center}
we can see that $t_k\geq2$ for any $k$ in all non-zero terms, so following the discussions above, we have
\begin{center}
$E(\sum\limits_{t=1}^N(\xi_t^2-\sigma^2))^{2q}\leq B_2 N^q$.
\end{center}
By far, we have proven that $(r_N(\lambda))^{2q}$ and $r_N^{2q}$ are $O_p(N^qm^{2q})$.
\begin{eqnarray*}
\left(Z_N(\lambda)\right)^{2q}&=&\left(\sum\limits_{t=1}^NZ_{N,t}(\lambda)\right)^{2q}\\
&=&\left(2\sum\limits_{t=1}^N\sum\limits_{v=1}^{m-1}\xi_t\xi_{t+v}w(vm^{-1})\cos(v\lambda)\right)^{2q}\\
&=&\sum\limits_{p=1}^{2q}\sum\limits_{t_1+\ldots+t_p=2q}\sum\limits_{j_1\cdots j_p}Z_{N,j_1}^{t_1}(\lambda)\cdots Z_{N,j_p}^{t_p}(\lambda).
\end{eqnarray*}
First, when $t_k\geq2$ for any $k$, $p\leq q$, the number of terms after expanding $Z_{N,j_1}^{t_1}$ is $(m-1)^{t_1}$. So for fixed $t_1,\ldots,t_p$, $Z_{N,j_1}^{t_1}\cdots Z_{N,j_p}^{t_p}$ contains no more than $(m-1)^{t_1+\cdots+t_p}=(m-1)^{2q}$ terms. What is more, for any fixed $p$, the number for all possible $j_1,\ldots,j_p$, of which the power satisfy $t_1+\cdots+t_p=2q$, is no more than $p!\,C_N^p$. So there are no more than $N^q m^{2q}$ terms, which means that $E\sum\limits_{p=1}^{q}\sum\limits_{t_1+\ldots+t_p=2q}\sum\limits_{j_1\cdots j_p}Z_{N,j_1}^{t_1}(\lambda)\cdots Z_{N,j_p}^{t_p}(\lambda)\leq B_1 N^qm^{2q}$.\par
When $p\leq q$ and some of $t_k$ are equal to 1, we assume that $t_{k_{s_1}}=\cdots=t_{k_{s_c}}=1$. For $t_{k_{s_1}}$, $j_{k_{s_1}}$ should satisfy $j_{k_{s_1}}-j_{k_{s_1}-1}<m$, that is, $\xi_{j_{k_{s_1}}}\in\{\xi_{j_{k_{s_1}-1}+1},\ldots,\xi_{j_{k_{s_1}-1}+m-1}\}$. If not, since $\xi_{j_{k_{s_1}}}$ is independent of $(Z_{N,j_1})^{t_1},\ldots,(Z_{N,j_{k_{s_1}-1}})^{t_{j_{k_{s_1}-1}}},\left(\xi_{j_{k_{s_1}}+1}+\cdots+\xi_{j_{k_{s_1}}+m-1}\right)$,\\
$\left(Z_{N,j_{k_{s_1}}+1}\right)^{t_{j_{k_{s_1}}+1}},\ldots ,\left(Z_{N,j_p}\right)^{t_p}$, we have
\begin{center}
$E(Z_{N,j_1})^{t_1}\cdots\left(Z_{N,j_p}\right)^{t_p}=0$.
\end{center}
So the number of all non-zero terms is no more than $C_{N}^{p-c} m^{c}$. For $(Z_{N,j_{k_{s_1}-1}})^{t_{j_{k_{s_1}-1}}}$, there are $(m-2)^{t_{j_{k_{s_1}-1}}}$ terms which contain $\xi_{j_{k_{s_1}}}$. So there are $(m-1)^{t_{j_{k_{s_1}-1}}}-(m-2)^{t_{j_{k_{s_1}-1}}}=O((m-2)^{t_{j_{k_{s_1}-1}}})$ terms which do not contain $\xi_{j_{k_{s_1}}}$. Then the number of non-zero terms is O($N^{p-c} m^c m^{t_1+\cdots+t_p-c})\leq O(N^q m^{2q})$. For $\xi_{j_{k_{s_1}}}$, if $Z_{N,j_{k_{1}}}^{t_{k_{1}}}, \ldots, Z_{N,j_{k_{d}}}^{t_{k_{d}}}$, contain $\xi_{j_{k_{s_1}}}$, then it is easy to see that there are totally no more than $m^{t_{k_1}+\cdots+t_{k_d}-1}$ terms which do not contain $\xi_{j_{k_{s_1}}}$. However, now for all $j_{k_{1}},\ldots,j_{k_{d}}$, $|j_{k_{v}}-j_{k_{s_1}}|\leq m-1$. So the number of all possible non-zero terms is $O(N^{p-c-d+1} m^{c+d} m^{t_1+\cdots+t_p-c})\leq O(N^qm^{2q})$. Hence $E(Z_{N,j_1})^{t_1}\cdots\left(Z_{N,j_p}\right)^{t_p}\leq B_1 N^qm^{2q}$. \par
When $p>q$, there are at least $p-q$ of $t_j=1$, so following the discussions above, $E(Z_{N,j_1})^{t_1}\cdots\left(Z_{N,j_p}\right)^{t_p}\leq B_1 N^qm^{2q}$. Therefore summarizing discussions above, $E(Z_N(\lambda))^{2q}=O_p(\frac{m^{2q}}{N^q})$.\par
\noindent Part 2, we prove that $\displaystyle \left|\frac{f(\lambda)-Ef_N(\lambda)}{f(\lambda)}\right|^{2q}=O(\frac{m^{2q}}{N^q})$.
\begin{eqnarray*}
\left|\frac{f(\lambda)-Ef_N(\lambda)}{f(\lambda)}\right|^{2q}&\leq& B\left|f(\lambda)-m\int W(m(\lambda-u))f(u)du\right|^{2q}\\
&&+B\left|m\int W(m(\lambda-u))(f(u)-EI_N(u))du\right|^{2q}
\end{eqnarray*}
And
\begin{eqnarray*}
\left|f(\lambda)-m\int W(m(\lambda-u))f(u)du\right|^{2q}&=&\left|\int (f(\lambda)-f(\lambda-um^{-1}))W(u)du\right|^{2q}\\
&\leq&\left(\int\left|f(\lambda)-f(\lambda-um^{-1})\right|W(u)du\right)^{2q}\\
&\leq& B_f m^{-2q}\left|\int |u|W(u)du\right|^{2q}
\end{eqnarray*}
Since $\alpha\geq \frac{1}{4}$, we have $\displaystyle \lim\limits_{N\rightarrow\infty}\frac{1}{m^{2q}}/\frac{m^{2q}}{N^{q}}<+\infty$.\par
And
\begin{eqnarray*}
\left|m\int W(m(\lambda-u))(f(u)-E I_N(u))du\right|^{2q}&=&\left|\int W(u)(f(\lambda-um^{-1})-E I_N(\lambda-um^{-1}))du\right|^{2q}\\
&\leq& \left(\int W(u)\left|f(\lambda-um^{-1})-E I_N(\lambda-um^{-1})\right|du\right)^{2q}
\end{eqnarray*}
Since by Woodroofe and Van Ness (1967),
\begin{center}
$\max_{|\lambda|\leq \pi}\left|f(\lambda)-E I_N(\lambda)\right|\leq B_{pe}\log N/N$,
\end{center}
where $B_{pe}$ is a constant. So
\begin{eqnarray*}
\left|m\int W(m(\lambda-u))(f(u)-E I_N(u))du\right|^{2q}&\leq& \left|\int W(u)\frac{B_{pe}\log N}{N}du\right|^{2q}\\
&\leq& B_{pe}^{2q}\frac{(\log N)^{2q}}{N^{2q}}.
\end{eqnarray*}
And $\displaystyle \lim\limits_{N\rightarrow\infty}\frac{(\log N)^{2q}}{N^{2q}}/\frac{m^{2q}}{N^{q}}\rightarrow0$. So we complete Part 2.\par
\noindent Part 3: we show that $E\left|\frac{f_N(\lambda)-Ef_N(\lambda)}{f(\lambda)}-( g_N(\lambda)-\sigma^2)\right|^{2q}=O_p(\frac{m^{2q}}{N^q})$.
\begin{eqnarray*}
&&\left|\frac{f_N(\lambda)-Ef_N(\lambda)}{f(\lambda)}-(g_N(\lambda)-\sigma^2)\right|^{2q}\\
&=&\left|\frac{f_N(\lambda)-E f_N(\lambda)-f(\lambda)(g_N(\lambda)-\sigma^2)}{f(\lambda)}\right|^{2q}\\
&\leq& B_f\left|f_N(\lambda)-E f_N(\lambda)-f(\lambda)(g_N(\lambda)-\sigma^2)\right|^{2q}\\
&=&B_f\left|m\int W(m(\lambda-u))\left(I_N(u)-E I_N(u)-f(u)(J_N(u)-E J_N(u))\right)du\right.\\
&&+\left.m\int W(m(\lambda-u))\left(f(\lambda)-f(u)\right)(J_N(u)-E J_N(u))du\right|^{2q}\\
&\leq&B_3\left|m\int W(m(\lambda-u))\left(I_N(u)-E I_N(u)-f(u)(J_N(u)-E J_N(u))\right)du\right|^{2q}\\
&&+B_3\left|m\int W(m(\lambda-u))\left(f(\lambda)-f(u)\right)(J_N(u)-E J_N(u))du\right|^{2q}\\
&=& R_1(\lambda)+R_2(\lambda),
\end{eqnarray*}
where $J_N(u)$ is the periodogram for $\xi_1,\ldots,\xi_N$, $B_3$ is a constant. After some manipulations (see Woodroofe and Van Ness 1967, Grenander and Rosenblatt 1957).
\begin{eqnarray*}
E R_1(\lambda)&=&E\frac{1}{N^{2q}}\left|\sum\limits_{r,s=-\infty}^{+\infty} a_ra_s d_{rs}(\lambda)\right|^{2q}\\
&\leq&\frac{1}{N^{2q}} \sum\limits_{p=1}^{2q}\sum\limits_{t_1+\cdots+t_p}\sum\limits_{\substack{r_1,s_1\\ \cdots \\ r_p,s_p}}E\left(\left|a_{r_1}a_{s_1}d_{r_1s_1}\right|^{t_1}\cdots\left|a_{r_p}a_{s_p}d_{r_ps_p}\right|^{t_p}\right)\\
&=&\frac{1}{N^{2q}}\sum\limits_{p=1}^{2q}\sum\limits_{t_1+\cdots+t_p}\sum\limits_{\substack{r_1,s_1\\ \cdots \\ r_p,s_p}}\left|a_{r_1}a_{s_1}\right|^{t_1}\cdots\left|a_{r_p}a_{s_p}\right|^{t_p}E\left(\left|d_{r_1s_1}\right|^{t_1}\cdots\left|d_{r_ps_p}\right|^{t_p}\right)\\
&\leq&\frac{1}{N^{2q}}\sum\limits_{p=1}^{2q}\sum\limits_{t_1+\cdots+t_p}\sum\limits_{\substack{r_1,s_1\\ \cdots \\ r_p,s_p}}\left|a_{r_1}a_{s_1}\right|^{t_1}\cdots\left|a_{r_p}a_{s_p}\right|^{t_p}\left(E\left|d_{r_1s_1}\right|^{2t_1}\cdots\left|d_{r_ps_p}\right|^{2t_p}\right)^{\frac{1}{2}}
\end{eqnarray*}
Now let us focus on $\left|d_{rs}\right|^{4q}$, that is $r_1=\cdots=r_p$, $s_1=\cdots=s_p$, and denote $R_\xi(v)$ as the autocovariance function of $\xi_j$.
\begin{eqnarray*}
\left|d_{rs}\right|^{4q}&=&\left|\sum\limits_{v_1=1-r}^{N-r}\sum\limits_{v_2=1-s}^{N-s} w((v_1-v_2+r-s)m^{-1})e^{-i(v_1-v_2+r-s)\lambda}(\xi_{v_1}\xi_{v_2}-R_\xi(v_1-v_2))\right.\\
& &\left.-\sum\limits_{v_1=1}^N\sum\limits_{v_2=1}^Nw((v_1-v_2+r-s)m^{-1})e^{-i(v_1-v_2+r-s)\lambda}(\xi_{v_1}\xi_{v_2}-R_\xi(v_1-v_2))\right|^{4q}\\
&=&\left|\sum\limits_{v_1=1-r}^{N-r}\sum\limits_{v_2=1-s}^{N-s}w((v_1-v_2+r-s)m^{-1})\cos((v_1-v_2+r-s)\lambda)(\xi_{v_1}\xi_{v_2}-R_\xi(v_1-v_2))\right.\\
&&+i\sum\limits_{v_1=1-r}^{N-r}\sum\limits_{v_2=1-s}^{N-s}w((v_1-v_2+r-s)m^{-1})\sin((v_1-v_2+r-s)\lambda)(\xi_{v_1}\xi_{v_2}-R_\xi(v_1-v_2))\\
& &-\sum\limits_{v_1=1}^N\sum\limits_{v_2=1}^Nw((v_1-v_2+r-s)m^{-1})\cos((v_1-v_2+r-s)\lambda)(\xi_{v_1}\xi_{v_2}-R_\xi(v_1-v_2))\\
& &\left.+i\sum\limits_{v_1=1}^N\sum\limits_{v_2=1}^Nw((v_1-v_2+r-s)m^{-1})\sin((v_1-v_2+r-s)\lambda)(\xi_{v_1}\xi_{v_2}-R_\xi(v_1-v_2))\right|^{4q}\\
&=&\left|Q_{rs}+iT_{rs}\right|^{4q}\leq 2^{2q}\left(Q_{rs}^{4q}+T_{rs}^{4q}\right)
\end{eqnarray*}
Again, we only need to prove that the expectation of $Q_{rs}^{4q}$ and $T_{rs}^{4q}$ are $O(N^{2q}m^{4q})$. And
\begin{eqnarray*}
E Q_{rs}^{4q}&=& E\left(\sum\limits_{v_1=1-r}^{N-r}\sum\limits_{v_2=1-s,v_2\neq v_1}^{N-s}w((v_1-v_2+r-s)m^{-1})\cos((v_1-v_2+r-s)\lambda)\xi_{v_1}\xi_{v_2}\right.\\
& &-\sum\limits_{v_1=1}^N\sum\limits_{v_2=1,v_2\neq v_1}^Nw((v_1-v_2+r-s)m^{-1})\cos((v_1-v_2+r-s)\lambda)\xi_{v_1}\xi_{v_2}\\
& &+\sum\limits_{v_1=\max(1-r,1-s)}^{\min(N-r,N-s)}\left(\xi_{v_1}^2-1\right)\left.-\sum\limits_{v_1=1}^{N}(\xi_{v_1}^2-1)\right)^{4q}\\
&\leq&2^{16q}\left(E\left(\sum\limits_{v_1=1-r}^{N-r}\sum\limits_{v_2=1-s,v_2\neq v_1}^{N-s}w((v_1-v_2+r-s)m^{-1})\cos((v_1-v_2+r-s)\lambda)\xi_{v_1}\xi_{v_2}\right)^{4q}\right.\\
&+&E\left(\sum\limits_{v_1=\max(1-r,1-s)}^{N-r,N-s}\left(\xi_{v_1}^2-\sigma^2\right)\right)^{4q}\\
&+&E\left(\sum\limits_{v_1=1}^N\sum\limits_{v_2=1,v_2\neq v_1}^Nw((v_1-v_2+r-s)m^{-1})\cos((v_1-v_2+r-s)\lambda)\xi_{v_1}\xi_{v_2}\right)^{4q}\\
&+&\left.E\left(\sum\limits_{v_1=1}^{N}(\xi_{v_1}^2-\sigma^2)\right)^{4q}\right)
\end{eqnarray*}
For those four terms in the last inequality above, following the proofs in Part 1, we can show that each of them is $O(N^{2q}m^{4q})$. Similarly to the discussions above, $E \left(T_{rs}(\lambda)\right)^{4q}=O(N^{2q}m^{4q})$. So we have $E\left|d_{rs}(\lambda)\right|^{4q}\leq O(N^{2q}m^{4q})$. When some of $(r_k,s_k)$ are different from each other, for example, $(r_1,s_1)\neq (r_k,s_k)$ for $k\geq2$, then in $E|d_{r_1s_1}(\lambda)|^{2t_1}\cdots|d_{r_ps_p}(\lambda)|^{2t_p}$, the number of non-zero terms should be less than $E\left|d_{rs}(\lambda)\right|^{4q}$ because some $\xi_j$ in $d_{r_1s_1}(\lambda)$ are not in $d_{r_ks_k}(\lambda)$. And now, when the power of $\xi_j$ is 1 in the expansion of $\left(d_{r_1s_1}(\lambda)\right)^{2t_1}$, we can not find any $\xi_j$ in $\left(d_{r_ks_k}(\lambda)\right)^{2t_k}$, so the expectations of these terms are 0. So,
\begin{eqnarray*}
E R_1(\lambda)\leq \frac{B N^{q}m^{2q}}{N^{2q}}\sum\limits_{p=1}^{2q}\sum\limits_{t_1+\cdots+t_p}\sum\limits_{\substack{r_1,s_1\\ \cdots \\ r_p,s_p}} \left|a_{r_1}a_{s_1}\right|^{t_1}\cdots\left|a_{r_p}a_{s_p}\right|^{t_p}
\end{eqnarray*}
Since $\sum_{j}|a_j|<+\infty$,
\begin{center}
$\sum\limits_{\substack{r_1,s_1\\ \cdots \\ r_p,s_p}}\left|a_{r_1}a_{s_1}\right|^{t_1}\cdots\left|a_{r_p}a_{s_p}\right|^{t_p}<+\infty$
\end{center}
for any given $t_1,\ldots,t_p$. Since $t_1+\cdots+t_p=2q$, $\displaystyle E R_1(\lambda)\leq B_1\frac{N^{q}m^{2q}}{N^{2q}}=B_1\frac{m^{2q}}{N^{q}}$.
\begin{eqnarray*}
R_2(\lambda)&=&B_2\left(\int W(u)\left(f(\lambda)-f(\lambda-um^{-1})\right)\left(J_N(\lambda-um^{-1})-E\left(J_N(\lambda-um^{-1})\right)\right)du\right)^{2q}\\
&=&B_3\left(\int W(u)^{1-\frac{1}{2q}}W(u)^{\frac{1}{2q}}\left(f(\lambda)-f(\lambda-um^{-1})\right)\times\right.\\
& &\left.\left(J_N(\lambda-um^{-1})-E\left(J_N(\lambda-um^{-1})\right)\right)du\right)^{2q}\\
&\leq&B_4\left(\int W(u)^{1-\frac{1}{2q}}W(u)^{\frac{1}{2q}}um^{-1}\left(J_N(\lambda-um^{-1})-E\left(J_N(\lambda-um^{-1})\right)\right)du\right)^{2q}\\
&\leq& B_4\left(\int \left(\left(W(u)\right)^{1-\frac{1}{2q}}\right)^{\frac{2q}{2q-1}}du\right)^{2q \frac{2q-1}{2q}}\times\\
& &\left(\int u^{2q}m^{-2q}W(u)\left(J_N(\lambda-um^{-1})-E J_N(\lambda-um^{-1})\right)^{2q}du\right)^{2q\frac{1}{2q}}\\
&=&B_4m^{-2q} \left|\int u^{2q}W(u)\left(J_N(\lambda-um^{-1})-E J_N(\lambda-um^{-1})\right)^{2q}du\right|\\
&=&B_4m^{-2q}\left|i^{-2q}\int  (iu)^{2q} W(u)\left(J_N(\lambda-um^{-1})-E J_N(\lambda-um^{-1})\right)^{2q}du\right|\\
&=&B_4m^{-2q}\left|\int \mathcal{F}(w^{(2q)})(u)\left(J_N(\lambda-um^{-1})-E J_N(\lambda-um^{-1})\right)^{2q}du\right|\\
&=&B_4m^{-2q}\int \mathcal{F}(w^{(2q)})(u)\left(\sum\limits_{t=1}^N\sum\limits_{v=1}^{N-t}\xi_t\xi_{t+v}e^{-i(\lambda-um^{-1})v}+\sum\limits_{t=1}^N(\xi_{t}^2-1)\right)^{2q}du
\end{eqnarray*}
where the second inequality holds because of H\"{o}lder inequality, and $\mathcal{F}(w^{(2q)})(u)$ is the Fourier transformation of the $2q$th derivative of $w(v)$. Since $w(v)=0$ for $v\geq 1$, we have $w^{(2q)}=0$ for $v\geq 1$, and
$w^{(2q)}(v)$ is bounded. So similar to Part 1, we have $\displaystyle R_2(\lambda)\leq B_5m^{-2q}\frac{m^{2q}}{N^q}=\frac{B_5}{N^q}$. Here we complete Part 3.\par
Since
\begin{eqnarray*}
\left|\frac{f_N(v)-f(v)}{f(v)}\right|^{2q}&\leq&\left|\frac{f_N(v)-E f_N(\lambda)+E f_N(\lambda)-f(\lambda)}{f(\lambda)}-(g_N(\lambda)-1)+(g_N(\lambda)-1)\right|^{2q}\\
&\leq&2^{6q}\left(\left|\frac{f_N(\lambda)-Ef_N(\lambda)}{f(\lambda)}-(g _N(\lambda)-\sigma^2)\right|^{2q}+\left|\frac{f(\lambda)-Ef_N(\lambda)}{f(\lambda)}\right|^{2q}\right.\\
& &\left.+\left|g_N(\lambda)-\sigma^2\right|^{2q}\vphantom{\left|\frac{f_N(\lambda)-Ef_N(\lambda)}{f(\lambda)}\right|}\right),
\end{eqnarray*}
we have
\begin{center}
$\displaystyle E\left|\frac{f_N(v)-f(v)}{f(v)}\right|^{2q}\leq B\frac{m^{2q}}{N^q}$.
\end{center}
Then, by Markov inequality,
\begin{eqnarray*}
P\left(\left|\frac{\hat{f}(\lambda)-f(\lambda)}{f(\lambda)}\right|\geq \epsilon \frac{m}{N^{1/2}}\right)&\leq&\frac{N^q}{m^{2q}\epsilon^{2q}}E\left|\frac{f_N(v)-f(v)}{f(v)}\right|^{2q}\\
&\leq&\frac{N^q}{m^{2q}\epsilon^{2q}}B\frac{m^{2q}}{N^q}=\frac{B}{\epsilon^{2q}}
\end{eqnarray*}
\begin{eqnarray*}
P\left(\max_{\lambda_j}\left|\frac{f_N(\lambda_j)-f(\lambda_j)}{f(\lambda_j)}\right|\geq\epsilon\frac{m}{n^{(q-1)/(2q)}}\right)&\leq&\sum\limits_{j=1}^NP\left(\left|\frac{f_N(\lambda_j)-f(\lambda_j)}{f(\lambda_j)}\right|\geq\epsilon\frac{m}{N^{(q-1)/(2q)}}\right)\\
&\leq&\sum\limits_{j=1}^N \frac{N^{q-1}}{\epsilon^{2q}m^{2q}}B_5\frac{m^{2q}}{N^q}=\frac{B_5}{\epsilon^{2q}}.
\end{eqnarray*}
So
\begin{eqnarray*}
P(\max\limits_{1\leq k< l\leq N}\max_{\lambda_i}\left|\frac{f_N(\lambda_i)-f(\lambda_i)}{f(\lambda_i)}\right|\geq\epsilon\frac{m}{N^{\frac{q-3}{2q}}})\leq\sum\limits_{k=1}^N\sum\limits_{l=1}^NB_5\frac{m^{2q}N^{q-3}}{N^{q-1}\epsilon^{2q}}=\frac{B_5m^{2q}}{\epsilon^{2q}}.
\end{eqnarray*}
Here we complete the proof of Lemma \ref{ch3l1}. \par
\begin{lemma}\label{ch3l2}
Suppose $\forall k$, $X_t=\sum\limits_{j=-\infty}^{+\infty} a_j(k)\epsilon_{t-j}$ satisfy Assumptions 1-8, then
\begin{equation}
f_0(\lambda)=\sum\limits_{k=1}^{K+1}\frac{N_k}{N} f_k(\lambda)\\
\end{equation}
\begin{equation}
\max\limits_{\lambda_i\in\Lambda}\left|\frac{\hat{f}_0(\lambda_i)-f_0(\lambda_i)}{f_0(\lambda_i)}\right|=O_p\left(\frac{m}{N^{(q-1)/2q}}\right)
\end{equation}
\end{lemma}
\noindent\textbf{Proof}:
\begin{eqnarray*}
\hat{f}_0(\lambda)&=&m\int_{-\infty}^{+\infty} W(m(u-\lambda))\hat{I}_0(u)du\\
&=&m\int_{-\infty}^{+\infty} W(m(u-\lambda))\frac{1}{N}\left|\sum\limits_{j=1}^N e^{-iju}X_j\right|^2du\\
&=&m\int_{-\infty}^{+\infty} W(m(u-\lambda))\frac{1}{N}\left|\sum\limits_{k=1}^{K+1}\sum\limits_{j=\tau_{k-1}^0+1}^{\tau_k^0}e^{-iju}X_j\right|^2du\\
&=&m\int_{-\infty}^{+\infty} W(m(u-\lambda))\frac{1}{N}\left(\sum\limits_{k=1}^{K+1}\left|\sum\limits_{j=\tau_{k-1}^0+1}^{\tau_k^0}e^{-iju}X_j\right|^2\right.\\
& &\left.+\sum\limits_{k_1\neq k_2}\sum\limits_{j=\tau_{k_1-1}^0+1}^{\tau_{k_1}^0}e^{-iju}X_j\sum\limits_{j=\tau_{k_2-1}^0+1}^{\tau_{k_2}^0}e^{iju}X_j\right)du\\
&=&m\int_{-\infty}^{+\infty} W(m(u-\lambda))\frac{1}{N}\left(\sum\limits_{k=1}^{K+1}\left|\sum\limits_{j=1}^{N_k}e^{-iju}X_j\right|^2\left|e^{i\tau_{k-1}^0u}\right|^2+\sum\limits_{k_1\neq k_2}\zeta_{k_1}(u)\bar{\zeta}_{k_2}(u)\right)du\\
&=&\sum\limits_{k=1}^{K+1}\hat{f}_{k}(\lambda)+\sum\limits_{k_1\neq k_2}\frac{m}{N}\int_{-\infty}^\infty W(m(u-\lambda))\zeta_{k_1}(u)\bar{\zeta}_{k_2}(u)du,
\end{eqnarray*}
where $\zeta_{k_1}(u)=\sum\limits_{j=\tau_{k_1-1}^0+1}^{\tau_{k_1}^0}e^{-iju}X_j$. So
\begin{eqnarray*}
& &P\left(\left|\frac{\hat{f}_0(\lambda)-f_0(\lambda)}{f_0(\lambda)}\right|\geq\epsilon\frac{m}{N^{1/2}}\right)\\
&\leq& P\left(B_f\left|\frac{\hat{f}_0(\lambda)-f_0(\lambda)}{f_0(\lambda)}\right|\geq\epsilon\frac{m}{N^{1/2}}\right)\\
&\leq& \sum\limits_{k=1}^{K+1}P\left(\frac{N_k}{N}\left|\hat{f}_k(\lambda)-f_k(\lambda)\right|\geq\epsilon\frac{m}{N^{1/2}(K+1)^2}\right)\\
& &+\sum\limits_{k_1\neq k_2}P\left(\left|\frac{m}{N}\int_{-\infty}^{+\infty} W(m(u-\lambda))\zeta_{k_1}(u)\bar{\zeta}_{k_2}(u)du\right|\geq \epsilon\frac{m}{N^{1/2}(K+1)^2}\right).
\end{eqnarray*}
By Lemma \ref{ch3l1}, we only need to show that
\begin{center}
$\displaystyle \left|\frac{m}{N}\int_{-\infty}^{+\infty} W(m(u-\lambda))\zeta_{k_1}(u)\bar{\zeta}_{k_2}(u)du\right|=O_p(\frac{m}{N^{1/2}})$.
\end{center}
What is more, assume $k_1<k_2$ without loss of generality, then we have
\begin{eqnarray*}
& &\left|\frac{m}{N}\int_{-\infty}^{+\infty} W(m(u-\lambda))\zeta_{k_1}(u)\bar{\zeta}_{k_2}(u)du\right|\\
&=&\frac{1}{N}\left|\int W(u)\left(\zeta_{k_1}(um^{-1}+\lambda)\bar{\zeta}_{k_2}(um^{-1}+\lambda)\right)du\right|\\
&=&\frac{1}{N}\left|\int W(u)\left(\sum\limits_{j=\tau_{k_1-1}^0}^{\tau_{k_1}^0}e^{-ij(um^{-1}+\lambda)}X_j\sum\limits_{j=\tau_{k_2-1}^0}^{\tau_{k_2}^0}e^{ij(um^{-1}+\lambda)}X_j\right)du\right|\\
&=&\frac{1}{N}\left|\sum\limits_{l=\tau_{k_2-1}^0+1-\tau_{k_1}^0}^{\tau_{k_2-1}+1-\tau_{k_1}^0+N_{k_1}}w(um^{-1})e^{il\lambda}\sum\limits_{j=\tau_{k_1}^0-(l-\tau_{k_2-1}^0+1-\tau_{k_1}^0)}^{\tau_{k_1}^0}X_jX_{j+l}\right.\\
& &+\sum\limits_{l=\tau_{k_2-1}+1-\tau_{k_1}+N_{k_1}+1}^{\tau_{k_2}-\tau_{k_1}}w(um^{-1})e^{il\lambda}\sum\limits_{j=\tau_{k_1-1}}^{\tau_{k_1}}X_jX_{j+l}\\
& &\left.+\sum\limits_{l=\tau_{k_2}-\tau_{k_1}+1}^{\tau_{k_2}-\tau_{k_1}-1}w(um^{-1})e^{il\lambda}\sum\limits_{j=\tau_{k_1-1}+1}^{\tau_{k_1}+1+N_{k_2}-N_{k_1}-(l-\tau_{k_2}-\tau_{k_1}+1)}X_jX_{j+l}\right|.
\end{eqnarray*}
Since $N_k>m$, $\forall k$, $w(u)=0$ for $|u|\geq 1$, so if $k_2-k_1>1$, $\zeta_{k_1}(u)\bar{\zeta}_{k_2}(u)=0$. So without loss of generality, we set $k_1=1$, $k_2=1$. Next, we separate the proofs into 2 parts, as in Lemma \ref{ch3l1}. \par
\noindent Part 1:
\begin{eqnarray*}
& &\frac{1}{N}\left|m\int W(m(u-\lambda))\zeta_1(u)\bar{\zeta}_2(u)\right|\\
&=&\frac{1}{N}\left|\int W(u)\sum\limits_{j=1}^{N_1}e^{-ij(um^{-1}+\lambda)}\xi_j\sum\limits_{j=N_1+1}^{N_2}e^{ij(um^{-1}+\lambda)}\xi_j\right|\\
&=&\frac{1}{N}\left|\sum\limits_{l=1}^{N_1}w(lm^{-1})e^{il\lambda}\sum\limits_{j=N_1-l+1}^{N_1-m}\xi_j\xi_{j+l}\right|.
\end{eqnarray*}
So
\begin{eqnarray*}
& &E\frac{1}{N^{2q}}\left|\sum\limits_{l=1}^{N_1}w(lm^{-1})e^{il\lambda}\sum\limits_{j=N_1-l+1}^{N_1-m}\xi_j\xi_{j+l}\right|^{2q}\\
&\leq&E\left(\frac{1}{N^{2q}}2^{2q}\left(\sum\limits_{l=1}^mw(lm^{-1})\cos(l\lambda)\sum\limits_{j=N_1-l+1}^{N_1-m}\xi_j\xi_{j+l}\right)^{2q}\right.\\
& &\left.+\frac{1}{N^{2q}}2^{2q}\left(\sum\limits_{l=1}^mw(lm^{-1})\sin(l\lambda)\sum\limits_{j=N_1-l+1}^{N_1-m}\xi_j\xi_{j+l}\right)^{2q}\right).\\
\end{eqnarray*}
So for two terms in the inequality above, following the proofs in Part 1 of Lemma \ref{ch3l1}, we have
\begin{center}
$\displaystyle E\left(\frac{1}{N^{2q}}2^{2q}\left(\sum\limits_{l=1}^mw(lm^{-1})\cos(l\lambda)\sum\limits_{j=n_1-l+1}^{n_1-m}\xi_j\xi_{j+l}\right)^{2q}\right)=O\left(\frac{m^{3q}}{N^{2q}}\right)$,\\
$\displaystyle E\left(\frac{1}{N^{2q}}2^{2q}\left(\sum\limits_{l=1}^mw(lm^{-1})\sin(l\lambda)\sum\limits_{j=n_1-l+1}^{n_1-m}\xi_j\xi_{j+l}\right)^{2q}\right)=O\left(\frac{m^{3q}}{N^{2q}}\right)$.
\end{center}
Part 2: Set $f_{11}(u)$, $f_{22}(u)$ satisfying $f_1(u)=f_{11}(u)\bar{f}_{11}(u)$, $f_2(u)=f_{22}(u)\bar{f}_{22}(u)$, where $\bar{f}_{22}$ denotes the conjugate of $f_{22}$. We show that $\displaystyle \left|\frac{m\int W(m(u-\lambda))\zeta_1(u)\bar{\zeta}_2(u)du}{f_{11}\bar{f}_{22}}- g_n(\lambda)\right|^{2q}=O_p(\frac{m^{2q}}{N^{q}})$.
\begin{eqnarray*}
& &\left|\frac{m\int W(m(u-\lambda))\zeta_1(u)\bar{\zeta}_2(u)du}{f_{11}\bar{f}_{22}}- g_n(\lambda)\right|^{2q}\\
&\leq&B_f\left|m\int W(m(u-\lambda))\zeta_1(u)\bar{\zeta}_2(u)du-f_{11}(\lambda)\bar{f}_{22}(\lambda)g_n(\lambda)\right|^{2q}\\
&\leq&\left|m\int W(m(u-\lambda))(\zeta_1(u)\bar{\zeta}_2(u)-f_{11}(u)\bar{f}_{22}(u))\right|^{2q}\\
& &+\left|m\int W(m(u-\lambda))(f_{11}(\lambda)\bar{f}_{22}(\lambda)-f_{11}(u)\bar{f}_{22}(u))J_1(u)J_2(u)du\right|^{2q}\\
&=&R_1(\lambda)+R_2(\lambda).
\end{eqnarray*}
Since $\forall k$, $f_k(\lambda)$ all satisfy uniform Lipschitz condition, then it is easy to see that $f_{11}(u)\bar{f}_{22}(u)$ also satisfies  uniform Lipschitz condition. Following the proofs in Part 3 of Lemma \ref{ch3l1}, we have
\begin{center}
$\displaystyle R_2(\lambda)=O_p(\frac{1}{N^{q}})$.
\end{center}
Following the proofs in Part 3 of Lemma \ref{ch3l1} again, we have $\displaystyle R_1(\lambda)=|\sum\limits_{r=-\infty,s=-\infty}^{+\infty} a_{r}(1)a_s(2)d_{rs}|^{2q}$, where
\begin{eqnarray*}
d_{rs}&=&\sum\limits_{v_1=1-r}^{\tau_1-r}\sum\limits_{v_2=\tau_1+1-s}^{\tau_2-s}\xi_{v_1}\xi_{v_2}e^{i(v_1-v_2)(um^{-1}+\lambda)}w((v_2-v_1)m^{-1})\\
& &-\sum\limits_{v_1=1}^{\tau_1}\sum\limits_{v_2=1}^{\tau_2}\xi_{v_1-r}\xi_{v_2-s}e^{i(v_2-v_1)(um^{-1}+\lambda)}w((v_2-v_1)m^{-1}).
\end{eqnarray*}
So following the proofs in Part 3 of Lemma \ref{ch3l1} again. we have
\begin{center}
$\displaystyle R_1(\lambda)=O_p(\frac{m^{2q}}{N^q})$.
\end{center}
Here we complete the proof of Lemma \ref{ch3l2}.\par
\begin{lemma}\label{ch3l3}
Assume Assumption 1-8, $\forall s=1,\cdots,K+1$, $\max\limits_{\tau_{s-1}^0\leq k< l\leq \tau_{s}^0}\vartheta_{kl}\sim O_p(N^{-\frac{q-3-2q\alpha}{2q}})$, when $l-k=N_{kl}\geq ml$. Here
\begin{center}
$\displaystyle \vartheta_{kl}=\frac{N_{kl}}{N}\int_{-\pi}^\pi \hat{f}_k^l(v)\log\left(\frac{st(\hat{f}_k)}{st(f_s)}\right)dv$
\end{center}
\end{lemma}
\noindent\textbf{Proof}: Set $m_{kl}=N_{kl}^\alpha$, and $B_1$ to $B_6$ are all constant. By Lemma \ref{ch3l1}, we have
\begin{center}
$\displaystyle P\left(\max\limits_{\lambda_i}\left|\frac{\hat{f}_{k}^l(\lambda_i)-f(\lambda_i)}{f(\lambda_i)}\right|^{2q}>\epsilon\frac{m^{2q}}{N^{(q-3)}}\right)\leq \frac{B_5}{N^{2}\epsilon^{2q}}$,
\end{center}
So denote $\displaystyle A_{kl}=\{\max\limits_{\lambda_i}\left|\frac{\hat{f}_{k}^l(\lambda_i)-f(\lambda_i)}{f(\lambda_i)}\right|\leq\epsilon\frac{m}{N^{(q-3)/(2q)}}\}$, then on set $A_{kl}$,
\begin{center}
$\displaystyle \left|\frac{\hat{f}_{k}^l(\lambda_i)-f(\lambda_i)}{f(\lambda_i)}\right|\leq\epsilon\frac{m}{N^{(q-3)/(2q)}}
\Rightarrow \max(0,1-\epsilon\frac{B_5m}{ N^{(q-3)/(2q)}})\leq \frac{\hat{f}_{kl}(\lambda_i)}{f_s(\lambda_i)}\leq 1+\epsilon\frac{B_5m}{ N^{(q-3)/(2q)}}$, $\forall k,l$, $\lambda_i$.
\end{center}
Since $\frac{\hat{f}_{kl}(\lambda_i)}{f_s(\lambda_i)}\geq0$,
\begin{center}
$\displaystyle \max(0,1-\epsilon\frac{B_5m}{ N^{(q-3)/(2q)}})\leq \frac{\hat{f}_{kl}(\lambda_i)}{f_s(\lambda_i)}\leq 1+\epsilon\frac{B_5m}{ N^{(q-3)/(2q)}}$.
\end{center}
So on set $A_{kl}$, we have
\begin{eqnarray*}
& &0\leq \frac{\hat{f}_{kl}(\lambda_i)}{f_s(\lambda_i)}\leq 1+\epsilon\frac{B_5m}{ N^{(q-3)/(2q)}}\\
&\Rightarrow& 0\leq \hat{f}_{kl}(\lambda_i)\leq f_s(\lambda_i)\left(1+\epsilon\frac{B_5m}{ N^{(q-3)/(2q)}}\right)\\
&\Rightarrow& 0\leq \hat{F}_{kl}(\pi)\leq F_s(\pi)\left(1+\epsilon\frac{B_6m}{ N^{(q-3)/(2q)}}\right).
\end{eqnarray*}
Then on $\bigcap\limits_{kl}A_{kl}$,
\begin{eqnarray*}
\max\limits_{kl}\vartheta_{kl}&=&\max\limits_{kl}\frac{N_{kl}}{N}\int_{-\pi}^\pi \hat{f}_k^l(u)\log\left(\frac{st(\hat{f}_k)}{st(f_s)}\right)du\\
&=&\max\limits_{kl}\frac{N_{kl}}{N}\int_{-\pi}^\pi \hat{f}_k^l(u)\log\left(\frac{F_s(\pi)}{\hat{F}_k(\pi)}\times\frac{\hat{f}_k(u)}{f_s(u)}\right)du\\
&\leq&\max\limits_{kl}\frac{N_{kl}}{N}\int_{-\pi}^\pi \max\limits_{u}\hat{f}_k^l(u)\log\left(\frac{F_s(\pi)}{\hat{F}_k(\pi)}\times\max\limits_{u}\frac{\hat{f}_k(u)}{f_s(u)}\right)du\\
&\leq&\max\limits_{kl}\frac{N_{kl}}{N}\int_{-\pi}^\pi \max\limits_{u}f_s(u)\left(1+\frac{B_5m\epsilon}{ N^{(q-3)/(2q)}}\right)\\
& &\log\left(\left(1+\frac{B_6m\epsilon}{ N^{(q-3)/(2q)}}\right)\left(1+\frac{B_5m\epsilon}{ N^{(q-3)/(2q)}}\right)\right)du\\
&\overset{\log(1+x)\leq x}{\leq}&\max\limits_{kl}\frac{N_{kl}}{N}\int_{-\pi}^\pi \max\limits_{u}f_s(u)\left(1+\frac{B_5m\epsilon}{ N^{(q-3)/(2q)}}\right)\left(\frac{B_6m\epsilon}{ N^{(q-3)/(2q)}}+\frac{B_5m\epsilon}{ N^{(q-3)/(2q)}}\right)du\\
&\leq& B_{7\epsilon}\frac{m}{N^{(q-3)/(2q)}}=B_{7\epsilon}N^{-\frac{q-3-2q\alpha}{2q}}.
\end{eqnarray*}
Here $B_{7\epsilon}$ is some constant containing $\epsilon$. So set $B_{8\epsilon}>B_{7\epsilon}$,
\begin{eqnarray*}
P\left(\max\limits_{kl}\vartheta_{kl}\geq B_{8\epsilon}N^{-\frac{q-3-2q\alpha}{2q}}\right)&=&P\left(\max\limits_{kl}\vartheta_{kl}\geq B_{8\epsilon} N^{-\frac{q-3-2q\alpha}{2q}}|\bigcap\limits_{kl}A_{kl}\right)P\left(\bigcap\limits_{kl}A_{kl}\right)\\
& &+P\left(\max\limits_{kl}\vartheta_{kl}\geq B_{8\epsilon}N^{\frac{q-3-2q\alpha}{2q}}|\bigcup\limits_{kl}A^c_{kl}\right)P\left(\bigcup\limits_{kl}A_{kl}^c\right)\\
&\leq&P\left(\max\limits_{kl}\vartheta_{kl}\geq B_{8\epsilon}N^{-\frac{q-3-2q\alpha}{2q}}|\bigcap\limits_{kl}A_{kl}\right)+P\left(\bigcup\limits_{kl}A_{kl}^c\right)\\
&\leq&0+\frac{B_5}{\epsilon^{2q}},
\end{eqnarray*}
and the last inequality can be arbitrarily small by choosing sufficiently large $\epsilon$. Here we complete proofs of Lemma \ref{ch3l3}.\par
\begin{lemma}\label{ch3l4}
Assume Assumption 1-8, $\forall s=1,\cdots,K+1$, $\displaystyle \max\limits_{1\leq k< l\leq N}\vartheta_{kl}\sim O_p(N^{-\frac{q-3-2q\alpha}{2q}})$, when $l-k=N_{kl}\geq ml$. Here
\begin{center}
$\displaystyle \vartheta_{kl}=\frac{N_{kl}}{N}\int_{-\pi}^\pi \hat{f}_k^l(v)\log\left(\frac{st(\hat{f}_k)}{st(f_s)}\right)dv$
\end{center}
\end{lemma}
\noindent\textbf{Proof}: Following the proofs of Lemma \ref{ch3l3}, this lemma can be easily obtained.\par
\vspace{1em}
\noindent\textbf{Proofs of Theorem \ref{esticonsis3}:}\\
Denote $\displaystyle A_{kl}=\left\{\omega: \max\limits_{\lambda_i}\left|\frac{\hat{f}_{k}^l(\lambda_i)-f(\lambda_i)}{f(\lambda_i)}\right|\leq\epsilon\frac{m}{N^{(q-3)/(2q)}},\forall j=1,\ldots,K\right\}$, where $\omega$ is the event in probability space $(\Omega, \mathcal{F}, P)$. Then $\forall \omega\in\bigcap\limits_{kl}A_{kl}$, we prove that $\hat{\kappa}_k\rightarrow \kappa_k^0$, $\forall j$. \par
For $\hat{\kappa}_k$, since $\hat{\kappa}_k(\omega)$ is bounded, $\forall k$, then there exists $\{n_s\}$ such that $\hat{\kappa}_{n_s}(\omega)\rightarrow\kappa_k^*$ on the subsequence. It follows from Lemma \ref{ch3l2} and \ref{ch3l3}, that
\begin{center}
$\displaystyle \frac{1}{N} R(\kappa_1^*,\ldots,\kappa_K^*)\leq\sum\limits_{i=1}^{K+1}(\kappa_i^*-\kappa_{i-1}^*)\int f_{\kappa_{i-1}^*\kappa_i^*}(v)\log \frac{st(f_{\kappa_{i-1}^*\kappa_i^*})}{st(f)}dv+O(N^{-\frac{q-3-2q\alpha}{2q}})$,
\end{center}
where $\lambda_0^*=0$, $\lambda_{K+1}^*=1$. If $\kappa_{i-1}\leq\kappa_{j-1}^*<\kappa_k<\cdots<\kappa_{i+k}<\kappa_j^*$, then by Lemma \ref{ch3l3}, we have
\begin{center}
$\displaystyle f_{\kappa_{j-1}^*\kappa_j^*}=\frac{\kappa_i-\kappa_{j-1}^*}{\kappa_j^*-\kappa_{j-1}^*}f_i+\frac{\kappa_{i+1}-\kappa_i}{\kappa_j^*-\kappa_{j-1}^*}f_{i+1}+\ldots+\frac{\kappa_j^*-\kappa_{j+k}}{\kappa_j^*-\kappa_{j-1}^*}f_{j+k+1}$\\
$\displaystyle F_{\kappa_{j-1}^*\kappa_j^*}=\frac{\kappa_i-\kappa_{j-1}^*}{\kappa_j^*-\kappa_{j-1}^*}F_i+\frac{\kappa_{i+1}-\kappa_i}{\kappa_j^*-\kappa_{j-1}^*}F_{i+1}+\ldots+\frac{\kappa_j^*-\kappa_{j+k}}{\kappa_j^*-\kappa_{j-1}^*}F_{j+k+1}$.
\end{center}
Since $\frac{1}{N}R(\kappa_1,\ldots,\kappa_K)=\sum\limits_{k=1}^{K+1}\int f_k(u)\log\frac{st(f_k)}{st(f)} du$.\par
So
\begin{eqnarray*}
& &\frac{1}{N}R(\lambda_{j-1}^*,\lambda_{j}^*)-\frac{1}{N}R(\lambda_{j-1}^*,\lambda_i)-\ldots-\frac{1}{N}R(\lambda_{i+k-1}^*,\lambda_{i+k})\\
&=&(\lambda_j^*-\lambda_i^0)\int f_i(u)\log\frac{st(f_{\lambda_{j-1}^*\lambda_j^*})}{st(f_i)}du+(\lambda_{i+1}^0-\lambda_i^0)\int f_{i+1}(u)\log\frac{st(f_{\lambda_{j-1}^*\lambda_j^*})}{st(f_{i+1})}du\\
& &+\cdots+(\lambda_j^*-\lambda_{i+k}^0)\int f_{i+k}(u)\log\frac{st(f_{\lambda_{j-1}^*\lambda_j^*})}{st(f_{i+k})}du+O(N^{-\frac{q-3-2q\alpha}{2q}})<0
\end{eqnarray*}
as $N\rightarrow\infty$, since every term above is always negative. If $\kappa_{j-1}< \kappa_k^*< \kappa_{k+1}^*<\kappa_j$, then
\begin{eqnarray*}
\frac{1}{N}R(\hat{\lambda}_{k-1},\hat{\lambda}_k)=(\lambda_j^*-\lambda_{j-1}^*)\int f_j(v)\log \frac{st(f_j)}{st(f)}dv+O(N^{-\frac{q-3-2q\alpha}{2q}}).
\end{eqnarray*}
So we have
\begin{eqnarray*}
\frac{1}{N}R(\hat{\kappa}_1,\ldots,\hat{\kappa}_K)=\frac{1}{N}R(\kappa^*_1,\ldots,\kappa_K^*)+O(N^{-\frac{q-3-2q\alpha}{2q}})-\frac{1}{N}R(\kappa_1^0,\ldots,\kappa_K^0)<0
\end{eqnarray*}
as $N\rightarrow\infty$.
This is a contradiction because $\hat{\kappa}_k$ are the maximizers of $R(\kappa_1,\ldots,\kappa_K)$. Since \begin{center}
$\displaystyle P\left(\bigcap\limits_{k,l}A_{kl}\right)=1-P\left(\bigcup\limits_{k,l}A^c_{kl}\right)\geq 1-\sum\limits_{l-k\geq ml}^NP\left(A_{kl}^c\right)=1-\frac{B_1}{\epsilon^{2q}}$.
\end{center}
so $\hat{\kappa}_j\overset{p}{\rightarrow} \kappa_j^0$.\par
\vspace{1em}
\noindent\textbf{Proofs of Theorem \ref{BICconsist3}:}\par
If $L<K$, then there should be a change point $\lambda_j^0$ that can not estimated consistently. Then following the proof in Theorem \ref{esticonsis3},
\begin{eqnarray*}
&&-\frac{1}{N}R(\hat{\kappa}_1,\ldots,\hat{\kappa}_L)+LC_N/N+\frac{1}{N}R(\hat{\kappa}_1,\ldots,\hat{\kappa}_K)-KC_N/N\\
&<&-c+O_p(N^{-\frac{q-3-2q\alpha}{2q}})-(K-L)C_N/N<0
\end{eqnarray*}
as $N\rightarrow\infty$. \par
If $\hat{K}>K$, then still every change point $\lambda_i$ should be estimated consistently. So, there should be a change point $\kappa_{j-1}<\kappa_k^*<\kappa_j$. Then by Lemma \ref{ch3l2} and \ref{ch3l3}, $\displaystyle R(\kappa_{j-1},\kappa_k^*,\kappa_j)-R(\kappa_{j-1},\kappa_{j})=O_p(N^{-\frac{q-3-2q\alpha}{2q}})$. So $BIC_L-BIC_K=O_p(N^{-\frac{q-3-2q\alpha}{2q}})+(K-L)C_N/N<0$, as $N\rightarrow\infty$. Here we complete the proofs of Theorem \ref{BICconsist3}.


\end{document}